\documentclass{cocv}
\usepackage{dsfont}
\usepackage{mathtools}
\usepackage{algorithmic}
\usepackage{algorithm}
\usepackage{braket,amsfonts}
\usepackage{mathtools}
\usepackage{enumitem}
\usepackage{enumitem}
\usepackage{array}
\usepackage{amssymb}
\usepackage{comment}
\newcommand{\R}{\mathds{R}}

\def\X{\rm{X}}

\def\Vad{{\rm V}^{\rm ad}} 
\def\VadStrict{\mathring{\rm{V}}^{\rm{ad}}}
\DeclareMathOperator*{\esssup}{ess\,sup}
\newcommand{\by}{\bar{y}}
\newcommand{\bq}{\bar{q}}
\newcommand{\bu}{\bar{u}}
\newcommand{\bx}{\bar{x}}
\newcommand{\bp}{\bar{p}}
\newcommand{\ue}{u_\epsilon}
\newcommand{\ze}{z_\epsilon}
\newcommand{\zd}{z_\delta}

\newcommand{\vd}{v_\delta}

\newcommand{\xdelta}{x_\delta}
\newcommand{\bue}{\bar{u}_\epsilon}
\newcommand{\bxe}{\bar{x}_\epsilon}
\newcommand{\bxen}{\bar{x}_{\epsilon_n}}
\newcommand{\bxenk}{\bar{x}_{\epsilon_{n_k}}}
\newcommand{\bpe}{\bar{p}_\epsilon}
\newcommand{\bpen}{\bar{p}_{\epsilon_n}}
\newcommand{\bpenk}{\bar{p}_{\epsilon_{n_k}}}
\newcommand{\buen}{\bar{u}_{\epsilon_n}}
\newcommand{\bues}{\bar{u}_{\epsilon_s}}
\newcommand{\buenk}{\bar{u}_{\epsilon_{n_k}}}

\newcommand{\pbe}{\bar{p}_\epsilon}

\newcommand\xs[1]{x^{1}}
\newcommand\supp{S}

\newcommand\norm[1]{\left\Vert#1\right\Vert}

\newcommand\const{\textrm{\textnormal{const}}}

\newcommand\st{\textrm{ s.t. }}

\def\xdt{{\rm{d}}t}
\def\xd{{\rm{d}}}
\def\xdn#1{{\rm{d}}#1}
\def\xD{{\rm{D}}}

\def\xBV{{\rm BV}}
\def\U{\textnormal{U}}
\def\AP{\textbf{AP} }
\def\meas{{\rm meas}}
\begin{document}
\title{Interior point methods in optimal control}
\author{Paul Malisani}\address{IFP Energies nouvelles, Applied Mathematics Department, 1 et 4 avenue de Bois-Pr\'eau, 92852 Rueil-Malmaison, France (\email{paul.malisani@ifpen.fr}).}

\date{24/04/2023}
\begin{abstract} This paper deals with Interior Point Methods (IPMs) for Optimal Control Problems (OCPs) with pure state and mixed constraints. This paper establishes a complete proof of convergence of IPMs for a general class of OCPs. Convergence results are proved for primal variables, namely state and control variables, and for dual variables, namely, the adjoint state, and the constraints multipliers. In addition, the presented convergence result does not rely on a strong convexity assumption. Finally, this paper compares the performances of a primal and a primal-dual implementation of IPMs in optimal control in three examples. \end{abstract}
\begin{resume} Cet article porte sur les méthodes de points-intérieurs pour des problèmes de commande optimale avec contraintes d'état pures et contraintes mixtes état-commande. Le résultat principal de cet article est une preuve de convergence de ces méthodes pour une classe générale de problèmes de commande optimale sous contraintes. Ce résultat de convergence est établi pour les variables primales, i.e., la commande et l'état et pour les variables duales, i.e. l'état-adjoint et les multiplicateurs des différentes contraintes. En outre, ce résultat de convergence ne repose pas sur une hypothèse de convexité forte du problème. Enfin, cet article compare les performances numériques entre une implantation primale et une implantation primale-duale des méthodes de points intérieurs en contrôle optimal autour de trois exemples. \end{resume}

\maketitle

\section{Introduction}
This paper deals with interior point methods (IPMs) for optimal control problems (OCPs) with pure state and mixed constraints. Numerous papers \cite{graichen,Weiser,Bonnans_log,lasdon67,maliOCAM,maliSIAM} deal with IPMs in optimal control.  However, unlike their counterpart in numerical optimization, they have yet to be fully theorized. In \cite{lasdon67,graichen}, the authors study interior point methods in optimal control with pure state constraints and prove the convergence of state and control variables using three strong assumptions. The first of these assumptions is the uniqueness of the optimal solution, the second is the strong convexity of the problem, and the last is the interiority of optimal trajectories concerning the constraints. In \cite{Bonnans_log}, the authors show the convergence of the adjoint state, the state, and the control for problems with control constraints of the form $a \leq u(t)\leq b$ and for nonlinear systems affine in the control variable. The proof of convergence in \cite{Bonnans_log} also relies on an assumption of uniqueness of the solution and an assumption of strong convexity of the problem. In \cite{Weiser}, the authors show the convergence of interior point methods in optimal control, with a primal-dual implementation, for problems with control constraints and using a strong Legendre-Clebsch condition, which, in this case, is equivalent to a strong convexity assumption \cite{BonnansHermant,FredericBonnans.2014}. In \cite{maliOCAM}, the authors establish sufficient conditions on the state and control penalties to ensure that the solutions of a penalized optimal control problem strictly satisfy the constraints and prove the convergence of state and control variables using an assumption of uniqueness of the optimal solution and a strong convexity assumption. In \cite{maliSIAM}, we recently proved IPMs convergence for primal variables (state and control), adjoint state, and multipliers associated with the pure-state and mixed constraints. This proof of convergence is established without using a strong convexity assumption but for OCPs with cost function, dynamics, and mixed constraints all affine in the control variable. Throughout this paper, we denote AP for Affine Property, an OCP where cost function, dynamics, and mixed constraints are affine in the control variable.\\
From a theoretical point of view, this paper generalizes IPMs convergence result from \cite{maliSIAM} to a general class of nonlinear OCPs without requiring additional assumptions and proves $\xLone$-convergence of the mixed-constraints multipliers versus weak $\ast$ convergence in \cite{maliSIAM}. From a numerical point of view, this paper provides detailed numerical examples using a primal implementation \cite{Bonnans_log,lasdon67,graichen,maliOCAM} and a primal-dual implementation as described in \cite{Weiser,maliSIAM}. These examples demonstrate interior-point methods' ease of use and versatility, particularly their primal-dual implementation.\\
The paper is organized as follows. \Cref{sec:problem_statement} contains the problem statement, the main assumptions, and the paper's main result. Some preliminary technical results are recalled in \cref{sec:preliminary_results}. In \cref{sec:existence_convergence}, we prove the existence of converging subsequences when using IPMs. In \cref{sec:proof_main_result}, we prove that the limit point of the aforementioned converging subsequences is a stationary point of the original optimal control problem, i.e., the limit point satisfies the first-order conditions of optimality \cite{Hartl,Maurer,Bonnans.2014}. In \cref{sec:primaldual}, we prove the convergence of the primal-dual implementation. Finally, in \cref{sec:algoExample}, we compare the performances of the primal and the primal-dual implementations on three challenging numerical examples. The first example is a state-constrained Van der Pol oscillator problem with quadratic cost, i.e., a cost which does not satisfy the \textbf{AP}. The second problem is a Zermelo's navigation problem, whose dynamics do not satisfy the \AP either. The third problem is the Goddard's problem \cite{seywald1993goddard}. This problem does satisfy the \textbf{AP} but is a classic numerical benchmark. The Python source codes are available in \cite{source_ifpen}.\\

\paragraph{\bf{Notations}:}
We denote $\R_-$ (resp. $\R_+$) the set of non-positive (resp. non-negative) real numbers. We denote $\mathds{N}_\ast$ (resp. $\R_\ast$) the set of non-zero natural integers (resp. real numbers). Given $p\in[1,+\infty]$, we denote $\xLn{p}(A;B)$ (or $\xLn{p}$) the Lebesgue spaces of functions from $A$ to $B$ and we denote $\norm{.}_{\xLn{p}}$ the corresponding $p$-norm. In addition, we also denote $\meas(.)$ the Lebesgue measure on $\R$. Given $p\in[1,+\infty]$, we denote $\xWn{{1,p}}(A;B)$ the Sobolev space of measurable functions from $A$ to $B$ with weak derivative in $\xLn{p}(A;B)$. Given $n\in[0,+\infty]$, we denote $\xCn{n}(A;B)$ (or $\xCn{n}$) the set of $n$-times continuously differentiable functions from $A$ to $B$. We denote $\xBV(A)$, the set of functions with bounded variations from $A$ to $\R$. We also denote $\mathcal{M}(A)$ the set of Radon measures on $A\subset\R$. The topological dual of a topological vector space $E$ is denoted $E^*$. Given a topological vector space $E$, we denote $\sigma(E,E^*)$ the weak topology on $E$ and $\sigma(E^*,E)$ the weak $\ast$ topology on $E^*$. Let $x_n,x\in E$, we denote $x_n\rightharpoonup x$ the weak convergence in $\sigma(E,E^*)$ and let $y_n,y\in E^*$, we denote $y_n\stackrel{\ast}{\rightharpoonup}y$ the weak $\ast$ convergence in $\sigma(E^*,E)$. For $x^*\in E^*$ and $x\in E$, we denote $\langle x^*,x\rangle$, the duality product. Given $f\in \xCn{{k\geq1}}(\R^n;\R)$ we denote $f'(.)$ the gradient of the function. Given $f\in \xCn{{k\geq1}}(\R^n \times \R^m;\R^p)$, we denote $f_x'(x,y):=\frac{\partial f}{\partial x}(x,y)\in \R^{p\times n}$ (resp.$f_y'(x,y):=\frac{\partial f}{\partial y}(x,y)\in \R^{p\times m}$) and we denote $f'_{i,x}:=\frac{\partial f_i}{\partial x}(x,y)$ (resp. $f'_{i,y}:=\frac{\partial f_i}{\partial y}(x,y)$). Given $f\in \xCn{{k\geq1}}(\R^n \times \R^m;\R)$, we denote $f'_{x,i}(x,y) := \left(f'_x(x,y)\right)_i$ (resp. $f'_{y,i}(x,y) := \left(f'_y(x,y)\right)_i$). We also denote $f''_{xy}(x,y):=\frac{\partial^2f}{\partial y \partial x}(x,y)$. Let $G:X\mapsto Y$ with $X,Y$ Banach spaces, we denote ${\rm D}G(x)$ the derivative of the mapping $G$ at point $x\in X$. The finite dimensional euclidean norm is denoted $\norm{.}$ and the scalar (resp. matrix) product between $x,y\in \R^n$ (resp. $x\in\R^{m\times n},y\in\R^n$) is denoted $x.y$. Given a set $E$, we denote $\vert E \vert$ its cardinal. We also denote $B_N(x,r)$ the closed ball of radius $r$ centered in $x$ for the topology induced by norm $N$. We denote $x[u,x^0]$ (or $x[u]$ if $x^0$ is fixed) the solution of the differential equations $\dot{x}=f(x,u)$ with initial condition $x^0$. Let $A,B$ be two spaces, we denote $A \cong B$ if $A$ and $B$ are isomorphic. Finally, we denote $\const(.)$ a positive finite constant depending on the parameters in argument.

 \newpage
\section{Problem statement and main result}
\label{sec:problem_statement}
\subsection{Optimal control problem}
The problem we are interested in consists in finding a solution  $(x,u)$ of the following Constrained Optimal Control Problem (COCP)
\begin{subequations}
\label{eq:all_orig_problem}
\begin{align}
\min_{(u,x)\in\U\times \X} J(x,u) & := \varphi(x(T)) + \int_0^T \ell(x(t), u(t))\xdt 
\label{eq:cost_orig}\\
    \dot{x}(t) & = f(x(t),u(t))  \label{eq:def_dynamics}\\
     0 & =h(x(0),x(T)) \label{eq:def_initial_condition}\\
    0 &\geq g(x(t)),\; \forall t\label{eq:def_state_const}\\
    0 &\geq c(x(t),u(t)),\;a.e. \label{eq:def_mixed_const}\\
    \U& := \xLinfty([0,T];\R^m)\label{eq:def_control_space}\\
    \X& := \xWn{{1,\infty}}([0,T];\R^n)\label{eq:def_state_space}
\end{align}
\end{subequations}
where the time horizon $T>0$ is fixed. In this setting, fixed final time is not a limitation since free end-time problems can be written as Problem \cref{eq:all_orig_problem} using the classical change in time variable $\tau:=t/T\in[0,1]$ and an augmented state whose dynamics writes $\begin{pmatrix}
  \dot{x}(\tau)&\dot{T}(\tau)
\end{pmatrix}^\top =\begin{pmatrix}
    T(\tau)f(x(\tau),u(\tau))& 0
\end{pmatrix}^\top$.
\begin{dfntn}
    \label{def:Vad}
    We denote $\Vad \subset \U \times \R^n$ the set of admissible controls and initial conditions as follows
    \begin{equation*}
    \label{eq:def_Vad}
        \Vad := \left\{(u,x^0)\in \U\times \R^n \textrm{ s.t. } \cref{eq:def_dynamics,eq:def_initial_condition,eq:def_mixed_const,eq:def_state_const} \textrm{ holds} \right\}
    \end{equation*}
    The set $\Vad$ is endowed with the following norm
    \begin{equation*}
        \label{eq:def_Vad_norm}
        \norm{(u,x^0)}_{\Vad} := \norm{u}_{\xLone} + \norm{x^0}
    \end{equation*}
    And, given $n\in \mathds{N}_\ast$, we denote $\VadStrict(n)$ the following set
    \begin{equation*}
        \label{eq:def_VadStrict}
        \VadStrict(n) := \left\{(u, x^0) \in \Vad \textrm{ s.t. } \begin{cases}
            g(x[u, x^0](t))&< 0,\; \forall t\\
            \esssup_t c(x[u, x^0](t), u(t))&\leq -\frac{1}{n}
        \end{cases}\right\}
    \end{equation*}
\end{dfntn}

\subsection{Main assumptions and technical definitions}
\begin{ssmptn}
\label{ass:C1_data}
The functions $\ell:\R^n\times\R^m\mapsto\R$, $f:\R^n\times\R^m\mapsto\R^n$, $g:\R^n\mapsto\R^{n_g}$, $c:\R^n\times \R^m \mapsto \R^{n_c}$, $h:\R^n\times \R^n\mapsto \R^{n_h}$ are at least twice continuously differentiable. In addition, $\ell$ and $c$ are convex with respect to their second argument, i.e., the control variable.
\end{ssmptn}

\begin{ssmptn}
    \label{ass:interior_accessibility}
    Any locally optimal solution $(x[\bu, \bx^0],\bu)$ such that $(\bu, \bx^0)\in\Vad$ satisfies the following interiority accessibility assumption
    \begin{equation*}
        \label{eq:def_interior_accessibility}
        (\bu, \bx^0) \in \Vad_\infty:=\textrm{\rm{cl}}_{\norm{.}_{\Vad}}\left(\liminf_{n} \VadStrict(n) \right)
    \end{equation*}
    where $\textrm{\rm{cl}}_{\norm{.}_{\Vad}}$ stands for the closure in the $\norm{.}_{\Vad}$-topology. 
\end{ssmptn}
\begin{rmrk}
  \Cref{ass:interior_accessibility} is essential for using IPMs in optimal control. If this assumption is not satisfied locally optimal solutions might not be accumulation points of strictly interior sequences and thus cannot be computed using IPMs. This assumption might not hold if, for example, the multipliers of Problem \cref{eq:all_orig_problem} are singular. Furthermore, we chose to endow $\Vad$ with the $\xLone$-topology in order to make this assumption as weak as possible.
\end{rmrk}

\begin{ssmptn}

\label{ass:bounded_if_state_set}
There exists $R_h < +\infty$ such that any locally optimal solution $(\bx,\bu)$ of Problem \cref{eq:all_orig_problem} satisfies
\begin{equation}
    \norm{\bx(0)},\norm{\bx(T)}\leq R_h
\end{equation}
\end{ssmptn}

\begin{ssmptn}
    \label{ass:bounded_control_set}
    There exists $R_u<+\infty$ such that for all $(u,x^0)\in \Vad$, we have
    \begin{equation}
        \label{eq:uniform_bounded_u}
        u \in B_{\xLinfty}(0,R_u)
    \end{equation}
    In addition, any solution of \cref{eq:def_dynamics} satisfies
    \begin{equation}
        \label{eq:ass_bibo}
        \forall R_v<+\infty,\exists R_x<+\infty,(u,x^0) \in B_{\norm{.}_{\Vad}}(0,R_v) \Rightarrow x[u,x^0] \in B_{\xLinfty}(0, R_x)
    \end{equation}
\end{ssmptn}

\begin{dfntn}[Sets of near state-saturated times and near-saturated indices]
For all $(u,x^0)\in \Vad$ from \cref{def:Vad}, for all $\delta \geq0$, and for all $n\in \mathds{N}_\ast$, we define the set of near state-saturated times (resp. mixed-saturated times), denoted $\supp^g_{u,x^0}$ (resp. $\supp^c_{u,x^0}$), as follows
\begin{align*}
    \supp^g_{u,x^0}(\delta) &:= \left\lbrace t\in [0,T] \st \max_i g_i(x[u, x^0](t))\geq -\delta\right\rbrace\\
    \supp^c_{u,x^0}(n) &:= \left\lbrace t\in [0,T] \st \max_i c_i(x[u, x^0](t), u(t))\geq -\frac{1}{n}\right\rbrace
\end{align*}
In addition, we define the set of near state-saturated indices (resp. mixed-saturated indices), denoted $I^g_{u,x^0}$ (resp. $I^c_{u,x^0}$) , as follows
\begin{align}
    I^g_{u,x^0}(t,\delta) &:= \left\lbrace i\in\{1,\dots,n_g\}\st g_i(x[u, x^0](t))\geq -\delta\right\rbrace\label{eq:def_indices_g}\\
    I^c_{u,x^0}(t,n) &:= \left\lbrace i\in\{1,\dots,n_c\} \st c_i(x[u, x^0](t), u(t))\geq -\frac{1}{n} \right\rbrace\label{eq:def_indices_c}
\end{align}
\end{dfntn}

\begin{ssmptn}
\label{ass:qualification_condition_mixed_constraint}
    For all $(u,x^0) \in \Vad$, the mixed constraints \cref{eq:def_mixed_const} satisfy the following qualification condition. There exists $\gamma>0$ and $n \in \mathds{N}_\ast$ such that
    \begin{equation*}
        \gamma \norm{\xi} \leq \norm{c'_{I^c_{u,x^0}(t,n),u}(x[u,x^0](t), u(t))^\top.\xi},\;\forall \xi \in \R^{\vert I^c_{u,x^0}(t,n)\vert},\;\textrm{a.a. } t\in [0,T]
    \end{equation*}
\end{ssmptn}

\begin{ssmptn}
    \label{ass:robinson_qualification}
    Let $(u, x^0)\in \Vad$ and denote $G:\U \times \R^n \mapsto \xCzero([0,T];\R^{n_g})$, $C:\U \times \R^n \mapsto \xLinfty([0,T];\R^{n_c})$ and $Z_{0T}: \U \times \R^n \mapsto \R^{n\times n}$ the constraints mapping defined by $G(v,y^0) := g(x[u + v,x^0 + y^0])$, $C(v, y^0):= c(x[u + v,x^0 + y^0], u+ v)$ and the linear mapping $Z_{0T}(v,y^0): = (z(0), z(T))$ where $z:[0,T]\mapsto \R^n$ is the unique solution of the following linearized state equation
    \begin{equation*}
        \begin{cases}
            \dot{z}(t) &= f'_x(x[u, x^0](t), u(t)). z(t) + f'_u(x[u, x^0](t), u(t)). v(t)\\
            z(0) &= y^0
        \end{cases}
    \end{equation*}
    There exists $\epsilon >0$ such that the following Robinson's constraints qualification conditions holds. 
    \begin{multline*}
        B_{ \xCzero \times \xLinfty \times \R^{n_h}}(0,\epsilon) \subset \begin{pmatrix}
              G(u, x^0)\\
              C(u,x^0) \\
              0
        \end{pmatrix} 
        + \begin{pmatrix}
              \xCzero([0,T],\R_-^{n_g})\\
              \xLinfty([0,T],\R_-^{n_c}) \\
              \{0\}
        \end{pmatrix}
        + \begin{pmatrix}
            \xD G(u, x^0)\\
            \xD C(u, x^0)\\
            \xD h(x[u, x^0](0),x[u, x^0](T))Z_{0T}
        \end{pmatrix} .(\U\times \R^n) 
    \end{multline*}
    As a consequence, the set of singular multipliers for Problem \cref{eq:all_orig_problem} is empty \cite[section 2.3.4., pp. 67]{Bonnans.Shapiro}.
\end{ssmptn}

\begin{rmrk}
Interested readers can refer to \cite{BonnansHermant} for sufficient conditions on pure-state and mixed constraints such that \cref{ass:robinson_qualification} holds. The mixed constraints qualification condition presented here is the same as in \cite{BonnansHermant}.
\end{rmrk}

\subsection{First-order necessary conditions of stationarity} This section recalls the first-order necessary conditions of optimality for Problem \cref{eq:all_orig_problem}. First, let us introduce the pre-Hamiltonian function.
\begin{dfntn}[pre-Hamiltonian]
\label{def:pre_hamiltonian}
Using \cref{ass:qualification_condition_mixed_constraint}, the set of singular multipliers is empty. Therefore, the pre-Hamiltonian $H:\R^n \times \R^m \times \R^n\mapsto \R$ of Problem \cref{eq:all_orig_problem} writes
\begin{equation*}
    \label{eq:def_pre_hamiltonian}
    H(x,u,p) := \ell(x,u) + p.f(x,u)
\end{equation*}
\end{dfntn}

\begin{dfntn}[Stationary point]
\label{def:pontryagin_extremal}
The trajectory $(\bx, \bu)$ with associated multipliers $(\bp, \bar{\mu}, \bar{\nu},\bar{\lambda})\in \xBV([0,T];\R)^n \times \xBV([0,T])^{n_g} \times  \xLinfty([0,T];\R^{n_c}_+) \times \R^{n_h}$, is a stationary point for Problem \eqref{eq:all_orig_problem} if it satisfies
\begin{subequations}
\label{eq:all_first_order}
\begin{align}
    \dot{\bar{x}}(t) =& f(\bar{x}(t),\bu(t))\label{eq:first_order_1}\\
    -\xdn{\bar{p}}(t) =& \begin{multlined}[t]\left[H_x'(\bx(t), \bu(t),\bp(t)) + \sum_{i=1}^{n_c}c_{i,x}'(\bx(t), \bu(t))\bar{\nu}_i(t)\right] \xdt +\sum_{i=1}^{n_g}g'_i(\bx(t))\xdn{\bar{\mu}}_i(t)\end{multlined}\label{eq:first_order_2}\\
    0=& H_u'(\bx(t),\bu(t), \bp(t))+\sum_{i=1}^{n_c} c_{i,u}'(\bx(t), \bu(t)) \bar{\nu}_i(t)\label{eq:first_order_3}\\
    0=&h(\bx(0), \bx(T)) \label{eq:first_order_4}\\
    0  =&\bp(0) + h_{x(0)}'(\bx(0),\bx(T))^\top.\bar{\lambda}\label{eq:first_order_5}\\
     0=&\bp(T) - \varphi'(\bx(T)) -  h_{x(T)}'(\bx(0),\bx(T))^\top.\bar{\lambda}\label{eq:first_order_6}\\
   0=& \int_0^T g_i(\bx(t))\xdn{\bar{\mu}}_i(t),\;\;i=1,\dots,n_g\label{eq:first_order_7}\\
    0=&\int_0^T c_i(\bx(t), \bu(t))\bar{\nu}_i(t)\xdt,\;\;i=1,\dots,n_c\label{eq:first_order_8}\\
    0 \leq & \xdn{\bar{\mu}}_i(t),\;\;i=1,\dots,n_g\label{eq:first_order_11}\\
    0 \leq &\bar{\nu}_i(t),\;\;i=1,\dots,n_c\label{eq:first_order_12}\\
    0 = &\bar{\mu}_i(T), \;\;,\;\;i=1,\dots,n_g\label{eq:first_order_13}
\end{align}
\end{subequations}
\end{dfntn}
Let $(\bx, \bu, \bp, \bar{\mu}, \bar{\nu}, \bar{\lambda})$ be a stationary point. Then, $( \bar{\mu}, \bar{\nu}, \bar{\lambda})$ is a generalized Lagrange multiplier (see \cite[Definition 2.7]{Bonnans.2014}). From \cref{ass:C1_data}, \cref{ass:interior_accessibility}, \cref{ass:bounded_if_state_set}, \cref{ass:bounded_control_set}, \cref{ass:qualification_condition_mixed_constraint}, \cref{ass:robinson_qualification} 
and using \cite[Theorem 3.1 and Corollary 3.2.]{Bonnans.2014}, we have that for any weak minimum (resp. Pontryagin minimum) $(\bx,\bu)$ of Problem \cref{eq:all_orig_problem}, there exists a generalized Lagrange multiplier (resp. generalized Pontryagin multiplier) $( \bar{\mu}, \bar{\nu}, \bar{\lambda})$ and their associated adjoint state $\bp$ such that $(\bx, \bu, \bp, \bar{\mu}, \bar{\nu}, \bar{\lambda})$ is a stationary point. Unfortunately, solving Problem \cref{eq:all_first_order} is a difficult task. Indeed, the dual variable $\xd\bar{\mu}$ associated with the state constraints appearing in \cref{eq:first_order_2,eq:first_order_7,eq:first_order_11,eq:first_order_13} is a Radon measure, therefore, in full generality, it can be decomposed in an absolutely continuous measure with respect to the Lebesgue measure, a discrete and finally a singular part. Computing this measure discrete and singular parts can be dramatically complex.
\subsection{Penalized Optimal Control Problem (POCP)}
To solve Problem \cref{eq:all_orig_problem}, we use an interior point method based on log-barrier functions defined as follows.
\begin{dfntn}[log-barrier function]
\label{def:log_barrier_function}
The log-barrier function $\psi:\R\mapsto\R$ is defined as follows
\begin{equation*}
\label{eq:def_log_barrier}
    \psi(x):=\begin{cases}
        -\log(-x)& \forall x<0\\
        +\infty&\textrm{otherwise}
    \end{cases}
\end{equation*}
\end{dfntn}
The penalized optimal control problem associated to \cref{eq:all_orig_problem} is defined as follows
\begin{subequations}
\label{eq:def_log_barrier_ocp}
\begin{align}
    \min_{(x,u)\in \X \times \U} J_\epsilon(x,u) &:=J(x,u) +\epsilon\int_0^T\left[\sum_{i=1}^{n_g}\psi\circ g_i(x(t)) +\sum_{i=1}^{n_c} \psi \circ c_i(x(t), u(t))\right] \xdt\\
    \dot{x}(t) & = f(x(t),u(t))\\
    0&=h(x(0),x(T))
\end{align}
\end{subequations}
The pre-Hamiltonian associated with this penalized problem is defined here after
\begin{dfntn}[Penalized pre-Hamiltonian]
\label{def:penalized_pre_hamiltonian}
The penalized pre-Hamiltonian $H^\psi:\R^n \times \R^m \times \R^n \times \R \mapsto \R$ of POCP \cref{eq:def_log_barrier_ocp} is defined by
\begin{equation*}
    \label{eq:def_penalized_pre_hamiltonian}
    H^\psi(x,u,p, \epsilon) := H(x,u,p)+ \epsilon\left(\sum_{i=1}^{n_g}\psi \circ g_i(x) +\sum_{i=1}^{n_c} \psi \circ c_i(x,u)\right)
\end{equation*}
\end{dfntn}

\begin{dfntn}[Penalized stationary point]
\label{def:penalized_pontryagin_extremal}
The trajectory $(\bxe, \bue)$ with associated multipliers $(\pbe, \bar{\lambda}_\epsilon)\in \xWn{{1,1}}([0,T];\R^n) \times \R^{n_h}$, is a penalized stationary point for Problem \eqref{eq:def_log_barrier_ocp} if it satisfies
\begin{subequations}
\label{eq:penalized_pontryagin_extremal}
\begin{align}
    \dot{\bar{x}}_{\epsilon}(t) =& f(\bxe(t),\bue(t))\label{eq:penalized_pontryagin_extremal_1}\\
    \dot{\bp}_\epsilon(t) =& - {H_x^{\psi}}'(\bxe(t),\bue(t),\pbe(t), \epsilon)\label{eq:penalized_pontryagin_extremal_2}\\
    0=&{H_u^{\psi}}'(\bxe(t),\bue(t),\pbe(t), \epsilon)\label{eq:penalized_pontryagin_extremal_3}\\
    0=&h(\bxe(0), \bxe(T))\label{eq:penalized_pontryagin_extremal_4}\\
    0 =& \bp_\epsilon(0) + h'_{x(0)}(\bxe(0), \bxe(T))^\top.\bar{\lambda}_{\epsilon}\label{eq:penalized_pontryagin_extremal_5}\\
    0 =& \bp_\epsilon(T) -\varphi'(\bxe(T))-  h'_{x(T)}(\bxe(0), \bxe(T))^\top.\bar{\lambda}_{\epsilon}\label{eq:penalized_pontryagin_extremal_6}
\end{align}
\end{subequations}
\end{dfntn}

\begin{ssmptn}
\label{ass:convex_hamiltonian}
Any penalized stationary point $(\bxe, \bue,\bpe, \bar{\lambda}_\epsilon)$ as defined in \cref{def:penalized_pontryagin_extremal} satisfies the following Hamiltonian convexity assumption
\begin{equation*}
    H''_{uu}(\bxe(t), \bue(t), \bpe(t)) \geq 0,\textrm{ a.e.}
\end{equation*}
\end{ssmptn}

\subsection{Contribution of the paper}
The main contribution of the paper is the generalization of the convergence result from \cite{maliSIAM} to general non linear optimal control problems. Throughout the rest of the paper and for the sake of readability, we will use the following notations for constraints multipliers.
\begin{align}
    \lambda_{\epsilon_n}^{g_i}& :=\epsilon_n \psi'\circ g_i(\bxen) \label{eq:def_lambda_g_mu}\\
    \lambda_{\epsilon_n}^{c_i} &:= \epsilon_n \psi'\circ c_i(\bxen, \buen)\label{eq:def_eta_c_nu}
\end{align}
We are now ready to state our main result.
\begin{thrm}
\label{thm:first_order_convergence}
Let $(\epsilon_n)_n$ be a sequence of decreasing positive parameters with $\epsilon_n\rightarrow0$. The associated sequence of penalized stationary points $(\bx_{\epsilon_n}, \bu_{\epsilon_n}, \bar{p}_{\epsilon_n}, \bar{\lambda}_{\epsilon_n})_n$ as defined in \cref{def:penalized_pontryagin_extremal} contains a subsequence converging to a stationary point $(\bx, \bu,\bp, \bar{\mu}, \bar{\nu},\bar{\lambda})$ of the original problem as defined in \cref{def:pontryagin_extremal}. The convergence is as follows
\begin{subequations}
\begin{align}
     \norm{\buenk - \bu}_{\xLone} & \rightarrow 0\label{eq:conv_u_L2}\\
      \norm{\bxenk -\bx}_{\xLinfty}&\rightarrow 0\label{eq:conv_x_Linfty}\\
      \norm{\bar{\lambda}_{\epsilon_{n_k}} - \bar{\lambda}}&\rightarrow 0 \label{eq:conv_lambda}\\
  \norm{\bpenk - \bar{p}}_{\xLone} &\rightarrow 0\label{eq:conv_p_L1}\\
\norm{ \lambda_{\epsilon_{n_k}}^{c_i} - \bar{\nu}_i}_{\xLone} & \rightarrow 0,\;\;\;i=1,\dots,n_c\label{eq:conv_weak_nu}\\
   \lambda_{\epsilon_{n_k}}^{g_i}\xdt&\stackrel{\ast}{\rightharpoonup} \xdn{\bar{\mu}}_i ,\;\;\;i=1,\dots,n_g\label{eq:weak_star_conv_pen_state}
\end{align}
\end{subequations}
\end{thrm}

\section{Preliminary results}
\label{sec:preliminary_results}
In this section, we gather the technical results established in \cite{maliSIAM} that do not rely on the $\AP$ of the OCP. These results can be used as is for OCP \cref{eq:all_orig_problem}.

\begin{prpstn}
\label{prop:major_Linf_L1}
Let $K<+\infty$, for all $(u_1, x^0_1), (u_2, x^0_2) \in B_{\norm{.}_{\Vad}}(0,K)$, there exists $\const(f) < +\infty $ such that
\begin{equation*}
    \parallel x[u_1, x_1^0]-x[u_2, x_2^0]\parallel_{\xLinfty} \leq \const(f)(\parallel u_1 - u_2\parallel_{\xLone} + \norm{x_1^0 - x_2^0})
\end{equation*}
\end{prpstn}
\begin{proof}
See \cite[Proposition 3.1]{maliSIAM}
\end{proof}

\begin{lmm}
\label{thm:interior_all}
There exists $K_g, K_c, K_p < +\infty$ such that for all $\epsilon \in(0, \epsilon_0]$, any associated locally optimal solution $(x[\ue, x^0_\epsilon], \ue)$ of Problem \cref{eq:def_log_barrier_ocp} and $(\bpe, \bar{\lambda}_\epsilon)$ their corresponding adjoint-state and initial-final constraint multiplier satisfy
\begin{align}
    \sup_{t\in [0,T]} g_i(x[\ue, x^0_\epsilon](t)) & < 0\label{eq:interiority_state}\\
    \norm{ \lambda_{\epsilon}^{g_i}}_{\xLone}& \leq K_g,\;\; i=1,\dots,n_g\label{eq:boundedness_lambda_g}\\
   \norm{ \lambda_{\epsilon}^{c_i}}_{\xLone}&\leq K_c,\;\; i=1,\dots,n_c\label{eq:boundedness_lambda_c}\\
   \norm{\bpe}_{\xLinfty} & \leq K_p\label{eq:boundedness_p}
\end{align}
\end{lmm}
\begin{proof}
We proved this result in \cite{maliSIAM} without using the \textbf{AP}. The proof is reproduced in \cref{sec:proof_lemma_siam} for this paper to be self-contained.
\end{proof}
\begin{lmm}
    \label{lem:conv_penalties}
    Let $(\bxen, \buen)_n$ be a sequence of local optimal solution of Problem \cref{eq:def_log_barrier_ocp}, then the following holds
    \begin{enumerate}[label=\roman*)]
        \item There exists $\bar{\mu}\in\xBV([0,T])^{n_g}$, with $\bar{\mu}(T)=0$ such that any sequence $(\lambda_{\epsilon_{n}}^{g_i})_n$ as defined in \cref{eq:def_lambda_g_mu} contains a subsequence satisfying \cref{eq:weak_star_conv_pen_state}.
        \item Any local optimal solution $(\bxen, \buen)$ of Problem \cref{eq:def_log_barrier_ocp} satisfies
        \begin{equation}
        \label{eq:interiority_mixed}
            \esssup_{t\in[0,T]} c_i(\bxe(t), \bue(t)) \leq - \epsilon/K_c
        \end{equation}
        \item Any element of the $(\lambda_{\epsilon_{n}}^{c})_n$ as defined in \cref{eq:def_eta_c_nu} satisfies
        \begin{equation}
        \label{eq:mixed_pen_in_Linfty}
            \lambda_{\epsilon_{n}}^{c} \in \xLinfty([0,T];\R^{n_c}_+)
        \end{equation}
        \item There exists $\bar{\nu}\in\xLinfty([0,T];\R^{n_c}_+)$ such that any sequence $(\lambda_{\epsilon_{n}}^{c})_n$ contains a subsequence satisfying 
        \begin{equation}
            \label{eq:conv_lambda_ci_nu_weakstar}
            \lim_{k \rightarrow + \infty} \lambda_{\epsilon_{n_k}}^{c_i} \stackrel{*}{\rightharpoonup} \bar{\nu}_i
        \end{equation}
    \end{enumerate}
\end{lmm}
\begin{proof}
The proof of this lemma is given in \cref{sec:add_proof} and is decomposed as follows. We proved $i)$ in \cite{maliSIAM} without using the \textbf{AP}, and the proof is reproduced in \cref{sec:proof_conv_mu}. The proof of $ii)$ (resp. $iii)$) is a straightforward adaptation of \cite{maliSIAM} established for \textbf{AP} problems and is given in \cref{sec:proof_of_2} (resp. \cref{sec:proof_of_3}) . Finally, $iv)$ is a direct consequence of \cref{eq:interiority_mixed} and of the weak $\ast$ compactness of $\xLinfty([0,T];\R^{n_c}_+)$.
\end{proof}

\section{Existence of converging sequences}
\label{sec:existence_convergence}
This section contains intermediate convergence results. We prove that the sequence of penalized stationary control $(\buen)_n$ contains a weakly converging subsequence to some $\bu$, that the associated state sequence $(\bxen)_n$ contains a uniformly converging subsequence, and that the associated sequence of adjoint state $(\bpen)_n$ contains an $\xLone$-converging subsequence. In addition, we will also prove that the sequences of the derivative of penalty functions also contain converging subsequences whose limit points belong to the dual space of their respective constraints. However, we call these convergence results \emph{intermediate} because when the problem does not satisfy the \textbf{AP}, there is no guarantee that any accumulation point of the sequence $(\bxen, \bpen)_n$ satisfies \cref{eq:first_order_1,eq:first_order_2}. An example illustrating this phenomenon is developed further in \cref{rmrk:example}. Nevertheless, these \emph{intermediate} convergence results are used in \cref{sec:proof_main_result} to prove the paper's  main result, namely, \cref{thm:first_order_convergence}.

\begin{lmm}
\label{prop:conv_weak_strong_un_xn}
There exists $(\by, \bu, \bq, \bar{\lambda}) \in \xCzero([0,T];\R^n)\times \xLtwo([0,T];\R^m) \times \\\xBV([0,T])^n  \times \R^{n_h}$ such that any sequence $(\bxen, \buen, \bpen, \bar{\lambda}_{\epsilon_n})_n$ of penalized stationary point as defined in \cref{def:penalized_pontryagin_extremal} for Problem \cref{eq:def_log_barrier_ocp} contains a subsequence satisfying
\begin{align}
    \lim_{k \rightarrow +\infty} \buenk &\rightharpoonup \bu \label{eq:weak_conv_un_bu} \\
    \lim_{k \rightarrow +\infty}\norm{\bx_{\epsilon_{n_k}} - \by}_{\xLinfty}&=0\label{eq:strong_xn_bx}\\
    \lim_{k \rightarrow +\infty} \norm{\bar{\lambda}_{\epsilon_{n_k}} - \bar{\lambda}}&=0\label{eq:weak_conv_lambda_bar}\\
    \lim_{k \rightarrow +\infty}\norm{\bpenk - \bq}_{\xLone}&=0\label{eq:strong_pn_bq}
\end{align} 
\end{lmm}
\begin{proof}
From \cref{ass:bounded_if_state_set} and \cref{ass:bounded_control_set}, the sequence $(\buen, \bxen^0)_n$ associated with $(\bxen, \buen)_n$ is $\xLtwo\times \R^n$ bounded. Therefore, there exists a weakly converging subsequence to $(\bu, \bx^0)\in \xLtwo\times \R^n$ which proves \cref{eq:weak_conv_un_bu}. To alleviate the notation we denote $(\buen, \bxen^0)_n$ the converging subsequence. Using this notation, we have
\begin{equation}
    x[\buen, \bxen^0](t_2) - x[\buen, \bxen^0](t_1) := \int_{t_1}^{t_2}f(x[\buen, \bxen^0](t), \buen(t)) \xdt
\end{equation}
First, $\forall t_1,t_2\in[0,T]$, and from Hölder inequality, we have
\begin{equation}
    \norm{x[\buen, \bxen^0](t_2)   - x[\buen, \bxen^0](t_1)}\leq\sup_n\norm{f(x[\buen, \bxen^0], \buen)}_{\xLtwo}\sqrt{\vert t_1 - t_2\vert} 
\end{equation}
Therefore, the sequence $(x[\buen, \bxen^0])_n$ is bounded and equicontinuous. From Arzela-Ascoli \cite[Theorem 1.3.8, p.33]{kurdila}, it contains a uniformly converging subsequence to some $\by$ which proves \cref{eq:strong_xn_bx}. In addition, \cref{eq:weak_conv_lambda_bar} is a direct consequence of the boundedness of the sequence $(\bar{\lambda}_{\epsilon_n})_n$. Now, let us prove \cref{eq:strong_pn_bq}. From \cref{thm:interior_all}, the sequence $({H^\psi_{x,i}}'(\bxen, \buen, \bpen, \epsilon_n))_n$ is uniformly $\xLone$-bounded. This sequence can be canonically injected into $\xCzero([0,T];\R^n)^{*}$ and is still uniformly bounded with respect to $\xCzero([0,T];\R^n)^{*}$ strong topology. From the weak $\ast$ compactness of the unit ball of $v$ (see \cite[Theorem 3.16]{brezis}) and since $\xCzero([0,T];\R^n)^{*}$ can be identified with the set of Radon measures $\mathcal{M}([0,T])^n$, there exists a subsequence and a measure $\xd \bar{q}\in \mathcal{M}([0,T])^n$ such that 
$$\lim_{k\rightarrow +\infty}{H^\psi_{x,i}}'(\bxenk, \buenk, \bpenk, \epsilon_{n_k}) \xd t\stackrel{\ast}{\rightharpoonup} \xd \bar{q}$$
Now, let us define $\bq\in \xBV([0,T])^n$ as follows
\begin{subequations}
\label{eq:def_bq}
\begin{align}
    \bq(T) - \bq(t) &:= \int^T_t \xd \bar{q}(s)\\
    \bq(T)& :=  \varphi'(\by(T)) +  h'_{x(T)}(\by(0), \by(T))^\top.\bar{\lambda}
    \end{align}
\end{subequations}
Then we have
\begin{multline*}
 \lim_{k\rightarrow +\infty }\norm{\bpenk(T) - \bpenk(t) - \bq(T) + \bq(t)} =\\
 \lim_{k\rightarrow +\infty } \norm{ \int^T_t {H^\psi_x}'(\bxenk(s), \buenk(s), \bpenk(s), \epsilon_{n_k})\xd s  - \xd \bar{q}(s) } = 0
\end{multline*}
Using \cref{eq:penalized_pontryagin_extremal_6,eq:strong_xn_bx,eq:weak_conv_lambda_bar,eq:def_bq}, we also have $\lim_{k\rightarrow+\infty} \norm{\bpenk(T)-\bq(T)}=0$, which yields $\lim_{k\rightarrow+\infty}\norm{ \bpenk(t)- \bq(t)}= 0$ and proves that $(\bpenk)_k$ pointwise converges to $\bq$. From \cref{thm:interior_all}, the sequence $(\bpenk)_k$ is uniformly $\xLinfty$-bounded which yields, from Lebesgue-Vitali's Theorem \cite[Theorem 4.5.4., pp. 268]{Bogachev.2007}
\begin{equation*}
     \lim_{k \rightarrow +\infty }\norm{\bpenk - \bq}_{\xLone} = 0
\end{equation*}
which proves \cref{eq:strong_pn_bq} and concludes the proof.
\end{proof}

\begin{rmrk}
    \label{rmrk:example}
    At this point, it is important to emphasize that since the problem does not satisfy the \textbf{AP}, the right-hand side of \cref{eq:penalized_pontryagin_extremal_1,eq:penalized_pontryagin_extremal_2} are, in full generality, non linear functions of $u$. This non linearity destroys weak convergence as illustrated in the following example. Let $\dot{x}(t) = u(t)^2$ and let $u_n(t):=\sin(nt)$. Then, we have
    \begin{align*}
        u_n&\rightharpoonup \bar{u}= 0\\
        (u_n)^2 &\rightharpoonup \frac{1}{2}\\
        x[u_n](t) - x[u_n](0)&\rightarrow \frac{t}{2}\\
        x[\bu](t) - x[\bu](0) &= 0
    \end{align*}
    Therefore $\lim_{n\rightarrow +\infty} x[u_n]\neq x[\bu]$. The proof of convergence for problems satisfying the  \textbf{AP} given in \cite{maliSIAM} extensively relies on the fact that for affine problems, if $u_n\rightharpoonup \bu$, then $\ell(x[u_n], u_n) \rightharpoonup \ell(x[\bu], \bu)$, $f(x[u_n], u_n) \rightharpoonup f(x[\bu], \bu)$ and $H'_x(x[u_n], u_n, p_n) \rightharpoonup H'_x(x[\bu], \bu, \bq)$ which is not true in the case we are interested in.
\end{rmrk}

\section{Proof of \cref{thm:first_order_convergence}}
To complete the proof of convergence it remains to prove:
\begin{itemize}
    \item The convergence of state and control variables, i.e., that \cref{eq:conv_u_L2,eq:conv_x_Linfty} hold. This is the object of \cref{sec:main_res_proof_x_u}
    \item The $\xLone$-convergence of the mixed constraints multipliers, i.e., that \cref{eq:conv_weak_nu} holds. This is the object of \cref{sec:main_res_proof_mixed_xLone}.
    \item The convergence of the adjoint state, i.e., that \cref{eq:conv_p_L1} holds. This is the object of \cref{sec:adjoint_state_conv_gen}.
    \item That the limit point of the sequences of locally penalized optimal solution satisfy the complementarity conditions, i.e., that \cref{eq:first_order_7,eq:first_order_8,eq:first_order_11,eq:first_order_12} hold. This the object of \cref{sec:main_res_conv_comp_gen}.
    \item That the limit point of the sequences of locally penalized optimal solution satisfy the Hamiltonian stationarity conditions, i.e., that \cref{eq:first_order_3} holds. This is the object of \cref{sec:proof_convergence_stationnarity}.
\end{itemize}

\label{sec:proof_main_result}
\subsection{Strong convergence of control and state variables}
\label{sec:main_res_proof_x_u}
Using \cref{ass:convex_hamiltonian} and from the convexity of $c(x,u)$ with respect to $u$ and since the penalty function \cref{eq:def_log_barrier} is strictly increasing, we have
\begin{equation*}
    \forall n>0, \textrm{ for a.e. } t\in[0,T], {H^{\psi}_{uu}}''(\bxen(t), \buen(t), \bpen(t), \epsilon_n) > 0
\end{equation*}
From the implicit function theorem \cite[Theorem 9.27, pp. 224-225]{Rudin}, for almost all time, there exists a mapping $\lambda_t$ such that
\begin{equation*}
    \buen(t) := \lambda_t(\bxen(t), \bpen(t), \epsilon_n)
\end{equation*}
By continuity of $\lambda_t$ and using \cref{prop:conv_weak_strong_un_xn} $(\bxen)_n$ (resp. $(\bpen)_n$) contains a pointwise converging subsequence, thus $(\buen)_n$ contains a subsequence pointwise converging to some $\hat{u} \in \xBV([0,T])^m$ satisfying
\begin{equation*}
    \hat{u}(t) = \lim_{k\rightarrow + \infty} \lambda_t(\bxenk(t), \bpenk(t), \epsilon_{n_k})=\lambda_t(\by(t), \bq(t), 0), \;\textrm{ a.e.}
\end{equation*}
Now, since there exists a control subsequence, denoted $(\bues)_s$, weakly converging to $\bu$, from Mazur's lemma \cite[lemma 10.19, pp. 350]{Renardy.2004}, there exists a function $S:\mathds{N} \mapsto \mathds{N}$ and a sequence of sets of real positive numbers $\left(\left\{\alpha[s]_k: k=s,\dots,S(s)\right\}\right)_s$ satisfying $\sum_{k=s}^{S(s)}\alpha[s]_k = 1$ and such that the sequence $(v_s)_s$ defined as follows
\begin{equation*}
    v_s := \sum_{k=s}^{S(s)}\alpha[s]_k \bar{u}_{\epsilon_k}
\end{equation*}
converges in $\xLtwo$-norm to $\bu$. Therefore, there exists a subsequence denoted $(v_m)_m$ converging almost everywhere to $\bu$. Now, for almost all $t\in[0,T]$, we have
\begin{equation*}
    \norm{\bu(t) - \hat{u}(t)} = \lim_{m\rightarrow +\infty} \norm{v_m(t) - \hat{u}(t)}
    \leq  \lim_{m\rightarrow +\infty} \sum_{k=m}^{S(m)}\alpha[m]_k\norm{\bu_{\epsilon_m}(t) - \hat{u}(t)}
    =0
\end{equation*}
which proves that $\hat{u}=\bu$ almost everywhere. As a consequence, there exists a subsequence $(\buenk)_k$ which converges almost everywhere to $\bu$. In addition, $(\buenk)_k$ being $\xLinfty$-bounded, using Lebesgue-Vitali's Theorem yields $\norm{\buenk - \bu}_{\xLone}\rightarrow 0$ wich proves \cref{eq:conv_u_L2}. In addition, using \cref{prop:major_Linf_L1} we have
\begin{equation*}
\label{eq:conv_bxenk_bx_proof}
    \norm{\bxenk - \bx}_{\xLinfty} \rightarrow 0
\end{equation*}
which proves \cref{eq:conv_x_Linfty}.

\subsection{$\xLone$-convergence of the mixed constraints multipliers}
\label{sec:main_res_proof_mixed_xLone}
From \cref{thm:interior_all} and \cref{lem:conv_penalties} the sequence $\bar{\lambda}_{\epsilon_{n}}^{c_i}$ is uniformly $\xLinfty$-bounded. The mapping $\psi'$ being differentiable, it is Lipschitz on bounded sets. Then, from the almost everywhere convergence of $\bxen$ and $\buen$ we get that there exists $\hat{\nu}_i$ such that
\begin{equation*}
    \bar{\lambda}_{\epsilon_{n}}^{c_i}\rightarrow \hat{\nu}_i,\;\;\textrm{ a.e.}
\end{equation*}
In addition, from the density of $\xLtwo([0,T];\R)$ in $\xLone([0,T];\R)$ and from \cref{lem:conv_penalties}, we have
\begin{equation*}
    \langle\bar{\lambda}_{\epsilon_{n}}^{c_i} - \bar{\nu}_i,f\rangle \rightarrow 0,\;\;\forall f\in \xLtwo([0,T];\R)
\end{equation*}
Thus, the sequence $\bar{\lambda}_{\epsilon_{n}}^{c_i} \rightharpoonup \bar{\nu}_i$ in the weak $\xLtwo$-topology. Therefore, using the same argument as in \cref{sec:main_res_proof_x_u} yields that $\hat{\nu}_i = \bar{\nu}_i$ and that there exists a subsequence $\bar{\lambda}_{\epsilon_{n_k}}^{c_i}$ satisfying
\begin{equation*}
    \lim_{k\rightarrow \infty}\norm{\bar{\lambda}_{\epsilon_{n_k}}^{c_i}- \bar{\nu}_i}_{\xLone} = 0
\end{equation*}
which proves \cref{eq:conv_weak_nu}.
\subsection{Convergence of adjoint state sequence}
\label{sec:adjoint_state_conv_gen}
Let us denote $q_n\in \xBV([0,T])^n $ the solution of 
\begin{align*}
    - \xd q_n(t)&= \begin{multlined}[t]\left[\ell_x'(\bxen(t),\buen(t))+f_x'(\bxen(t),\buen(t)).q_n(t)+\sum_{i=1}^{n_c}c_{i,x}'(\bxen(t),\buen(t))\bar{\nu}_i(t)\right]\xdt\\
+\sum_{i=1}^{n_g}g_i'(\bxen(t))\xd \bar{\mu}_i(t)\end{multlined}\\
 q_n(0)&=-h'_{x(0)}(\bxen(0), \bxen(T))^\top.\bar{\lambda} \\
 q_n(T)&=\varphi'(\bxen(T))+h'_{x(T)}(\bxen(0), \bxen(T)) ^\top .\bar{\lambda}
\end{align*}
The functions $\ell,f,g,c, h$ and $\varphi$ being $\xCtwo$ and $(\bxen, \buen)_n$ being uniformly bounded, the mapping $(\buen, \bxen(0))\mapsto q_n$ is a continuous mapping from $\xLone([0,T];\R^m)\times \R^n \mapsto \xBV([0,T])^n$, hence, from the strong $\xLone$-convergence (resp. $\xLinfty$-convergence) of $(\buen)_n$ (resp. $(\bxen)_n$), we have
\begin{equation}
    \lim_{n\rightarrow +\infty} \norm{q_n-\bp}_{\xBV} = 0
    \label{eq:conv_pq_BV}
\end{equation}
In addition, since $q_n(T) = \bpen(T)$, we have
\begin{align*}
    \norm{q_n(t)-\bpen(t)}\leq&\begin{multlined}[t] \Bigg{\|}\int^T_t f_x'(\bxen(s),\buen(s)).(q_n(s)-\bpen(s))\xd s + \sum_{i=1}^{n_g}\int_t^Tg_i'(\bxen(s))\left(\xd \bar{\mu}_i(s) - \lambda^{g_i}_{\epsilon_n}(s)\xd s\right)\\
    + \sum_{i=1}^{n_c}\int_t^Tc_{i,x}'(\bxen(s), \buen(s))\left(\bar{\nu}_i(s)- \lambda_{\epsilon_n}^{c_i}(s)\right)\xd s\Bigg{\|}
    \end{multlined}\\
    \leq&\begin{multlined}[t]\const(f)\int^T_t \norm{q_n(s)-\bpen(s)}\xd s + \sum_{i=1}^{n_g}\norm{\int_t^Tg_i'(\bxen(s))\left(\xd \bar{\mu}_i(s)-\lambda^{g_i}_{\epsilon_n}(s)\xd s\right)}\\
+\sum_{i=1}^{n_c}\norm{\int_t^Tc_{i,x}'(\bxen(s), \buen(s))\left[\bar{\nu}_i(s)-\lambda_{\epsilon_n}^{c_i}(s)\right]\xd s}
    \end{multlined}
\end{align*}
Now, let us define $h_n\in \xLone([0,T];\R_+)$ as follows
\begin{equation*}
h_n(t):=\sum_{i=1}^{n_g}\norm{\int_t^Tg_i'(\bxen(s))\left(\xd\mu_i(s)-\lambda^{g_i}_{\epsilon_n}(s)\xd s\right)}
    +\sum_{i=1}^{n_c}\norm{\int_t^Tc_{i,x}'(\bxen(s), \buen(s))\left[\bar{\nu}_i(s)-\lambda_{\epsilon_n}^{c_i}(s)\right]\xd s}
\end{equation*}
thus $\norm{q_n(t)-\bpen(t)}\leq\const(f)\int^T_t \norm{q_n(s)-\bpen(s)}\xd s + h_n(t)$. From Grönwall inequality \cite[Lemma A.1, p.651]{Khalil}, we have $\norm{q_n(t)-\bpen(t)}\leq \const(f,T)\int^T_th_n(s)\xd s$. Using the weak $*$ convergence of $\xd \bar{\mu}$, the  $\xLone$-convergence of $\bar{\nu}$, the $\xLinfty$-convergence of $\bxen$ and the $\xLone$-convergence of $\buen$ yields the pointwise convergence of $h_n\rightarrow 0$. In addition, from the boundedness of $(h_n)_n$ there exists a subsequence such that
\begin{equation*}
  \lim_{k\rightarrow +\infty}\norm{q_{n_k}(t)-\bpenk(t)}\leq \const(f,T)\int^T_t\lim_{k\rightarrow +\infty}h_{n_k}(s)\xd s=0
\end{equation*}
$q_{n_k}$ pointwise converges to $\bpenk$ and since both are bounded, from Lebesgue-Vitali's Theorem, we have
$$\lim_{k\rightarrow +\infty}\norm{q_{n_k}-\bpenk}_{\xLone} = 0$$ 
Gathering with \cref{eq:conv_pq_BV} yields 
\begin{equation*}
\label{eq:conv_penk_bq_proof}
\norm{\bp -\bpenk}_{\xLone} \rightarrow 0
\end{equation*}
which proves \cref{eq:conv_p_L1}.

\subsection{Complementarity conditions satisfaction for the limit point of the sequence of penalized stationary points}
\label{sec:main_res_conv_comp_gen}
Up to this point, we have proved convergence of the control, the state, the adjoint state and the constraints multipliers. Let us prove that this limit point satisfies the complementarity conditions from the first-order conditions of optimality \cref{eq:first_order_7,eq:first_order_8,eq:first_order_11,eq:first_order_12}.
\subsubsection{Complementarity conditions for state constraints}
Now, let us prove that $\bar{\mu}$ satisfies conditions \cref{eq:first_order_7,eq:first_order_11}. From \cref{thm:interior_all}, we have $\bar{\lambda}_{\epsilon_{n_k}}^{g_i}(t)>0$, $\forall t\in[0,T]$ and $\forall \epsilon_{n_k}>0$. Therefore, $\forall \phi\in \xCzero([0,T];\R_+)$ one has
\begin{equation*}
\int_0^T \phi(t) \xd \bar{\mu}_i(t) =\lim_{k\rightarrow +\infty}\int_0^T  \phi(t) \bar{\lambda}_{\epsilon_{n_k}}^{g_i}(t) \xdt \geq 0
\end{equation*}
which proves that $\bar{\mu}$ satisfies the non negativity condition \cref{eq:first_order_11}. Finally, let us prove that $\bar{\mu}$ satisfies the complementarity condition \cref{eq:first_order_7}. From \cref{eq:def_lambda_g_mu}, we have $g_i(\bxenk(t))\bar{\lambda}_{\epsilon_{n_k}}^{g_i}(t) = - \epsilon_{n_k}$ hence
\begin{equation}
\label{eq:conv_gin_lambdan_to_zero}
    \lim_{k\rightarrow+\infty}  \int_0^Tg_i(\bxenk(t))\bar{\lambda}_{\epsilon_{n_k}}^{g_i}(t)\xdt = \lim_{k\rightarrow+\infty}-\epsilon_{n_k} T=0
\end{equation}
From the continuity of $g_i$, the sequence $(g_i(\bxenk))_k$ uniformly converges to $g_i(\bx)$. In addition, from \cref{thm:interior_all}, the sequence $\bar{\lambda}_{\epsilon_{n_k}}^{g_i}(t)$ is uniformly $\xLone$-bounded, hence
\begin{equation}
\label{eq:conv_gistar_gicsue}
   \lim_{k\rightarrow +\infty}\left\vert\int_0^T\left(g_i(\bx(t))-g_i(\bxenk(t))\right)\bar{\lambda}_{\epsilon_{n_k}}^{g_i}(t) \xdt \right\vert \leq
   \lim_{k\rightarrow +\infty}\norm{g_i(\bx)-g_i(\bxenk)}_{\xLinfty}\norm{\bar{\lambda}_{\epsilon_{n_k}}^{g_i}}_{\xLone}=0
\end{equation}
Gathering \cref{eq:conv_gin_lambdan_to_zero,eq:conv_gistar_gicsue} yields
\begin{equation*}
    \lim_{k\rightarrow +\infty}\int_0^Tg_i(\bx(t))\bar{\lambda}_{\epsilon_{n_k}}^{g_i}(t)\xdt=  \lim_{k\rightarrow +\infty}\int_0^Tg_i(\bxenk(t))\bar{\lambda}_{\epsilon_{n_k}}^{g_i}(t)\xdt = 0
\end{equation*}
which in turns gives
\begin{equation*}
    \int_0^T g_i(\bx(t))\xd \bar{\mu}_i(t)  = \lim_{k\rightarrow +\infty}\int_0^Tg_i(\bx(t))\bar{\lambda}_{\epsilon_{n_k}}^{g_i}(t) \xdt
    =\lim_{k\rightarrow +\infty}\int_0^Tg_i(\bxenk(t))\bar{\lambda}_{\epsilon_{n_k}}^{g_i}(t) \xdt
    =0
\end{equation*}
and proves that $\bar{\mu}$ satisfies the complementarity condition \cref{eq:first_order_7}.
\subsubsection{Complementarity conditions for mixed constraints}
Now, let us prove that $\bar{\nu}$ satisfies conditions \cref{eq:first_order_8,eq:first_order_12}. From \cref{lem:conv_penalties} we have $\bar{\lambda}_{\epsilon_{n_k}}^{c_i}(t)>0$, $\forall t\in[0,T]$ and $\forall \epsilon_{n_k}>0$ which proves that $\bar{\nu}$ satisfies the non negativity condition \cref{eq:first_order_12}. Finally, let us prove that $\bar{\nu}$ satisfies the complementarity condition \cref{eq:first_order_8}. From \cref{eq:def_eta_c_nu}, we have $c_i(\bxenk(t), \buenk(t))\bar{\lambda}_{\epsilon_{n_k}}^{c_i}(t) = - \epsilon_{n_k}$ hence
\begin{equation}
\label{eq:conv_cin_etan_to_zero}
    \lim_{k\rightarrow+\infty}  \int_0^Tc_i(\bxenk(t), \buenk(t))\bar{\lambda}_{\epsilon_{n_k}}^{c_i}(t)\xdt = \lim_{k\rightarrow+\infty}-\epsilon_{n_k} T=0
\end{equation}
From the Lipshitz-continuity of $c_i$ over compact sets, the sequence $(c_i(\bxenk, \buenk))_k$  $\xLone$-converges to $c_i(\bx, \bu)$. In addition, from \cref{thm:interior_all} and \cref{lem:conv_penalties}, the sequence $\bar{\lambda}_{\epsilon_{n_k}}^{c_i}$ is uniformly $\xLinfty$-bounded, hence
\begin{multline}
\label{eq:conv_cistar_cicsue}
   \lim_{k\rightarrow +\infty}\left\vert\int_0^T\left(c_i(\bx(t), \bu(t))-c_i(\bxenk(t), \buenk(t))\right)\bar{\lambda}_{\epsilon_{n_k}}^{c_i}(t) \xdt \right\vert \leq \\
   \lim_{k\rightarrow +\infty}\norm{c_i(\bx, \bu)-c_i(\bxenk, \buenk)}_{\xLone}\norm{\bar{\lambda}_{\epsilon_{n_k}}^{c_i}}_{\xLinfty}=0
\end{multline}
Gathering \cref{eq:conv_cin_etan_to_zero,eq:conv_cistar_cicsue} yields
\begin{equation*}
\label{eq:conv_cista_lambda_ci}
    \lim_{k\rightarrow +\infty}\int_0^Tc_i(\bx(t), \bu(t))\bar{\lambda}_{\epsilon_{n_k}}^{c_i}(t)\xdt=  \lim_{k\rightarrow +\infty}\int_0^Tc_i(\bxenk(t), \buenk(t))\bar{\lambda}_{\epsilon_{n_k}}^{c_i}(t)\xdt = 0
\end{equation*}
which in turns gives
\begin{equation*}
    \int_0^T c_i(\bx(t), \bu(t)) \bar{\nu}_i(t)\xdt  = \lim_{k\rightarrow +\infty}\int_0^Tc_i(\bx(t), \bu(t))\bar{\lambda}_{\epsilon_{n_k}}^{c_i}(t) \xdt
    =\lim_{k\rightarrow +\infty}\int_0^Tc_i(\bxenk(t), \buenk(t))\bar{\lambda}_{\epsilon_{n_k}}^{c_i}(t) \xdt
    =0
\end{equation*}
which proves that $\bar{\nu}$ satisfies the complementarity condition \cref{eq:first_order_8}.

\subsection{Stationarity condition of $H'_u$ for the limit points of the sequence of penalized stationary points}
\label{sec:proof_convergence_stationnarity}
Finally, the last step of the proof of \cref{thm:first_order_convergence} consists in proving that the limit point of the sequence of penalized stationary points satisfies the Hamiltonian stationarity conditions \cref{eq:first_order_3}.
\begin{multline}
\label{eq:last_eq_proof}
    \lim_{n\rightarrow +\infty }\norm{{H_u^{\psi}}'(\bxen, \buen, \bpen, \epsilon_n) - H'_u(\bx, \bu, \bp) - c'_u(\bx, \bu).\bar{\nu}}_{\xLone} \leq \\
    \lim_{n\rightarrow +\infty }\bigg[\norm{\ell'_u(\bxen, \buen)-\ell'_u(\bx, \bu)}_{\xLone} +  \norm{f'_u(\bxen, \buen)^\top . \bpen -f'_u(\bx, \bu)^\top . \bp }_{\xLone} \\
    +\sum_{i=1}^{n_c}\norm{\bar{\lambda}_{\epsilon_{n}}^{c_i} c_{i,u}'(\bxen, \buen) - \bar{\nu}_i c_{i,u}'(\bx, \bu)}_{\xLone}\bigg]
\end{multline}
Since all the terms in the right-hand side of \cref{eq:last_eq_proof} converges to zero, using \cref{eq:penalized_pontryagin_extremal_3} proves \cref{eq:first_order_3} and concludes the proof.
\section{Primal-dual convergence}
\label{sec:primaldual}
\begin{crllr}
    \label{thm:primal_dual}
    Let $(\epsilon_n)_n$ be a sequence of decreasing positive parameters with $\epsilon_n\rightarrow0$ and let
    \begin{multline}
       (\bxen, \buen, \bpen, \bar{\lambda}_{\epsilon_n}^{g}, \bar{\lambda}_{\epsilon_n}^c,\bar{\lambda}_{\epsilon_n})_n\in\\
       \xWn{{1,\infty}}([0,T];\R^n) \times \textnormal{U} \times \xWn{{1,1}}([0,T];\R^n) \times \xLone([0,T];\R^{n_g}_+) \times \xLinfty([0,T];\R^{n_c}_+) \times \R^{n_h}
    \end{multline}
be a solution of the following primal-dual TPBVP
\begin{subequations}
\label{eq:all_primal_dual}
\begin{align}
    \dot{\bar{x}}_{\epsilon_n}(t)=& f(\bxen(t),\buen(t))\label{eq:primal_dual_1}\\
    \dot{\bar{p}}_{\epsilon_n} =&
    -H'_x(\bxen(t), \buen(t), \bpen(t))-g'(\bxen(t))^\top .\bar{\lambda}_{\epsilon_n}^{g}(t) - c_{x}'(\bxen(t), \buen(t))^\top .\bar{\lambda}_{\epsilon_n}^{c}(t)\label{eq:primal_dual_2}\\
    0=&H'_u(\bxen(t), \buen(t), \bpen(t)) + c_{u}'(\bxen(t), \buen(t))^\top .\bar{\lambda}_{\epsilon_n}^{c}(t) \label{eq:primal_dual_3}\\
    0=&{\rm FB}(\bar{\lambda}_{\epsilon_n}^{g_i}(t), g_i(\bxen(t)), \epsilon_n)\label{eq:primal_dual_4}\\
    0=&{\rm FB}(\bar{\lambda}_{\epsilon_n}^{c_i}(t),c_i(\bxen(t), \buen(t)),\epsilon_n)\label{eq:primal_dual_5}\\
    0 =&h(\bxen(0),\bxen(T))\label{eq:primal_dual_6}\\
    0=& \bpen(0)  +h'_{x(0)}(\bxen(0), \bxen(T))^\top.\bar{\lambda}_{\epsilon_n}\label{eq:primal_dual_7}\\
    0 =&\bpen(T) -\varphi'(\bxen(T))- h'_{x(T)}(\bxen(0), \bxen(T))^\top.\bar{\lambda}_{\epsilon_n}
\end{align}
\end{subequations}
where ${\rm FB}:\R \times \R, \times \R_+ \mapsto \R_+$ is the Fisher-Burmeister complementarity function defined as follows
\begin{equation}
    {\rm FB}(x, y, \epsilon) := x - y - \sqrt{x^2 + y^2 + 2 \epsilon}
\end{equation}
Then $(\bxen, \buen, \bpen, \bar{\lambda}_{\epsilon_n}^{g}, \bar{\lambda}_{\epsilon_n}^{c}, \bar{\lambda}_{\epsilon_n})_n$ contains a subsequence satisfying \cref{eq:conv_u_L2,eq:conv_x_Linfty,eq:conv_lambda,eq:conv_p_L1,eq:conv_weak_nu,eq:weak_star_conv_pen_state}

\end{crllr}

\begin{proof}
From \cref{thm:interior_all}, we have $g_i(\bxen(t))<0$ and  $c_i(\bxen(t),\buen(t))<0$ for all $\epsilon_n>0$. Therefore \cref{eq:primal_dual_4} is equivalent to $\bar{\lambda}_{\epsilon_n,i}^{g}(t) =-\epsilon_n/g_i(\bxen(t))$ and \cref{eq:primal_dual_4} is equivalent to $\bar{\lambda}_{\epsilon_n,i}^{c}(t) =-\epsilon_n/c_i(\bxen(t), \buen(t))$. Combining with \cref{eq:primal_dual_2,eq:primal_dual_3} proves that any solution of \cref{eq:all_primal_dual} is also solution of \cref{eq:penalized_pontryagin_extremal} and using \cref{thm:first_order_convergence} concludes the proof.
\end{proof}

\newpage
\section{Examples}
\label{sec:algoExample}

This section will solve three challenging optimal control problems using a primal and a primal-dual implementation of interior-point methods in optimal control. The primal algorithm is the same as described in \cite{Bonnans_log,lasdon67,graichen,maliOCAM}. The primal-dual algorithm is the same as in \cite{Weiser,maliSIAM}. Both algorithms are recalled in \cref{sec:primal_dual_algos}. In this section, for the sake of readability and since there is no possible confusion between function space and the image space of these functions, we will not explicitly write the time dependency of the variables.

\subsection{Primal and primal-dual algorithms}
\label{sec:primal_dual_algos}
In the following, for all $\epsilon >0$,  we denote
\begin{equation}
    {\rm S_P}(\epsilon):=(\bxe, \bpe, \bue,\bar{\lambda}_\epsilon)
\end{equation}
any solution of \cref{eq:penalized_pontryagin_extremal}. Using this notation, the classical primal algorithm is presented in \cref{alg:primal}.

\begin{algorithm}[H]
\caption{Primal algorithm for optimal control problems}
\label{alg:primal}
\begin{algorithmic}[1]
\STATE{Define $\epsilon_0>0,\; \alpha \in(0,1),\;\textrm{tol}=o(1),\;k=0 $}
\STATE{Initialize $(x,p,u,\lambda)$}
\STATE{Compute ${\rm S_P}(\epsilon_{0})$ from initialization variables} 
\WHILE{$\epsilon_k > \textrm{tol}$}
\STATE{Compute ${\rm S_P}(\epsilon_{k+1})$ using ${\rm S_P}(\epsilon_k)$ as initialization} 
\STATE{$\epsilon_{k+1}\gets \alpha\epsilon_k$}
\STATE{$k\gets k+1$}
\ENDWHILE
\RETURN ${\rm S_P}(\epsilon_k)$
\end{algorithmic}
\end{algorithm}
\noindent
In addition, for all $\epsilon >0$,  we denote
\begin{equation}
    {\rm S_{PD}}(\epsilon):=(\bxe, \bue, \bpe, \bar{\lambda}_{\epsilon}^{g}, \bar{\lambda}_{\epsilon}^{c},\bar{\lambda}_{\epsilon}) 
\end{equation}
any solution of \cref{eq:all_primal_dual}. The corresponding primal-dual algorithm as presented in \cite{Weiser,maliSIAM} is recalled in \cref{alg:primal_dual}.
\begin{algorithm}[H]
\caption{Primal-dual algorithm for optimal control problems}
\label{alg:primal_dual}
\begin{algorithmic}[1]
\STATE{Define $\epsilon_0>0,\; \alpha \in(0,1),\;\textrm{tol}=o(1),\;k=0 $}
\STATE{Initialize $(x,p,u,\lambda^g, \lambda^c, \lambda)$}
\STATE{Compute ${\rm S_{PD}}(\epsilon_{0})$ from initialization variables} 
\WHILE{$\epsilon_k > \textrm{tol}$}
\STATE{Compute ${\rm S_{PD}}(\epsilon_{k+1})$ using ${\rm S_{PD}}(\epsilon_k)$ as initialization} 
\STATE{$\epsilon_{k+1}\gets \alpha\epsilon_k$}
\STATE{$k\gets k+1$}
\ENDWHILE
\RETURN ${\rm S_{PD}}(\epsilon_k)$
\end{algorithmic}
\end{algorithm}

\subsection{Parameterization of the decay ratio $\alpha$}
To make a fair comparison between \cref{alg:primal} and \cref{alg:primal_dual}, the decay ratio $\alpha$ appearing in both methods must be chosen carefully. Setting the same value for $\alpha$ in both methods prevents comparing the capacity of both algorithms to handle small decay ratios. Therefore, for each example and each algorithm, the decay ratio is set at the minimal value that allows the algorithm to converge.

\subsection{Van der Pol oscillator}
\subsubsection{Problem presentation and two-point boundary value problems}
The first example is the Van der Pol oscillator which consists in solving the following problem
\begin{equation*}
    \min_{(x,u)}\left[J(x,u):=  \int_0^4 x_1^2 + x_2 ^2 + u^2 \xdt\right]
\end{equation*}
under the following constraints
\begin{align*}
    \dot{x}_1 &= x_2\\
    \dot{x_2} & = - x_1 + x_2\left(1 - x_1^2\right) + u\\
    x_1(0), x_2(0) & = 1\\
    x_1(4)^2 + x_2(4)^2 & = 0.04\\
    -0.4 - x_1 &\leq 0\\
    u-1 &\leq 0\\
    -u &\leq 0
\end{align*}
The cost function is quadratic with respect to the control variable, therefore the \textbf{AP} does not hold and IPMs convergence for this problem stems from \cref{thm:first_order_convergence}. One can verify that any primal solution ${\rm S_P}(\epsilon)$ is solution of the following TPBVP
\begin{align*}
        \dot{x}_1 & = x_2\\
        \dot{x}_2 & =  - x_1 + x_2\left(1 - x_1^2\right) + u\\
        \dot{p}_1 & = - 2 x_1 + p_2 + 2 p_2 x_2 x_1\\
        \dot{p}_2 & =  - 2 x_2 - p_1(t) + p_2 (x_1^2 - 1) - \frac{\epsilon}{x_2+0.4}\\
        0 & = 2 u + p_2 + \frac{\epsilon}{1 -u} - \frac{\epsilon}{1 + u}\\
        x_1(0), x_2(0) & = 1\\
    x_1(4)^2 + x_2(4)^2 & = 0.04\\
    p_1(4) x_2(4) - p_2(4) x_1(4) & = 0
\end{align*}
and that any primal-dual solution ${\rm S_{PD}}(\epsilon)$ is solution of 
\begin{align*}
     \dot{x}_1 & = x_2\\
        \dot{x}_2 & =  - x_1 + x_2\left(1 - x_1^2\right) + u\\
        \dot{p}_1 & = - 2 x_1 + p_2 + 2 p_2 x_2 x_1\\
        \dot{p}_2 & =- 2 x_2 - p_1 + p_2 (x_1^2 - 1) + \lambda^g\\
        0 & = 2 u + p_2 + \lambda^c_1 - \lambda^c_2\\
        0 & = {\rm FB}(\lambda^g, -0.4 - x_1, \epsilon)\\
        0 & = {\rm FB}(\lambda^c_1, u - 1, \epsilon)\\
        0 & ={\rm FB}(\lambda^c_2, -1 - u, \epsilon)\\
        x_1(0), x_2(0) & = 1\\
    x_1(4)^2 + x_2(4)^2 & = 0.04\\
    p_1(4) x_2(4) - p_2(4) x_1(4) & = 0
\end{align*}

\subsubsection{Initialization}
For both algorithms, the initial values of time, state, adjoint state, control, and the penalty parameter are as follows.
\begin{align*}
    {\rm time}&\gets \begin{bmatrix} 0, &0.1, &\dots, &3.9, &4\end{bmatrix}\\
    x_1, x_2&\gets\begin{bmatrix} 1, &\dots, &1\end{bmatrix} \\
    p_1, p_2, u&\gets \begin{bmatrix} 0, &\dots, &0\end{bmatrix}\\
    \epsilon_0 & \gets 1
\end{align*}
And the constraints multipliers for the primal-dual algorithm are initialized as follows.
\begin{equation*}
    \lambda^g,\lambda^c_1,\lambda^c_2 \gets \begin{bmatrix}0, &\dots, &0\end{bmatrix}
\end{equation*}
\subsubsection{Results}
\begin{figure}[H]
    \centering
    \includegraphics[width=\textwidth]{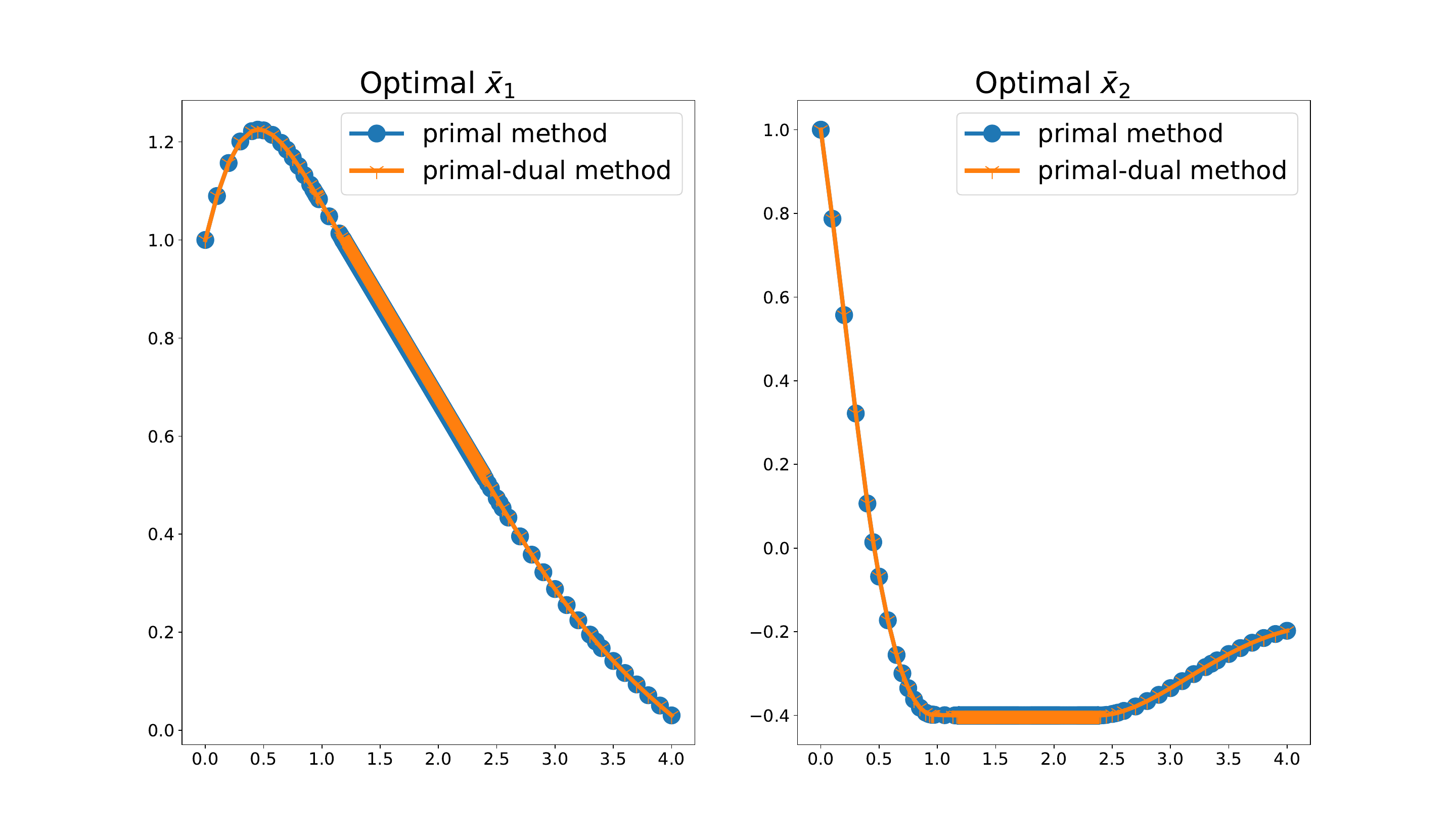}
    \caption{Optimal states $x_1$ and $x_2$ using primal and primal-dual methods}
    \label{fig:vdp_state}
\end{figure}

\begin{figure}[H]
    \centering
    \includegraphics[width=\textwidth]{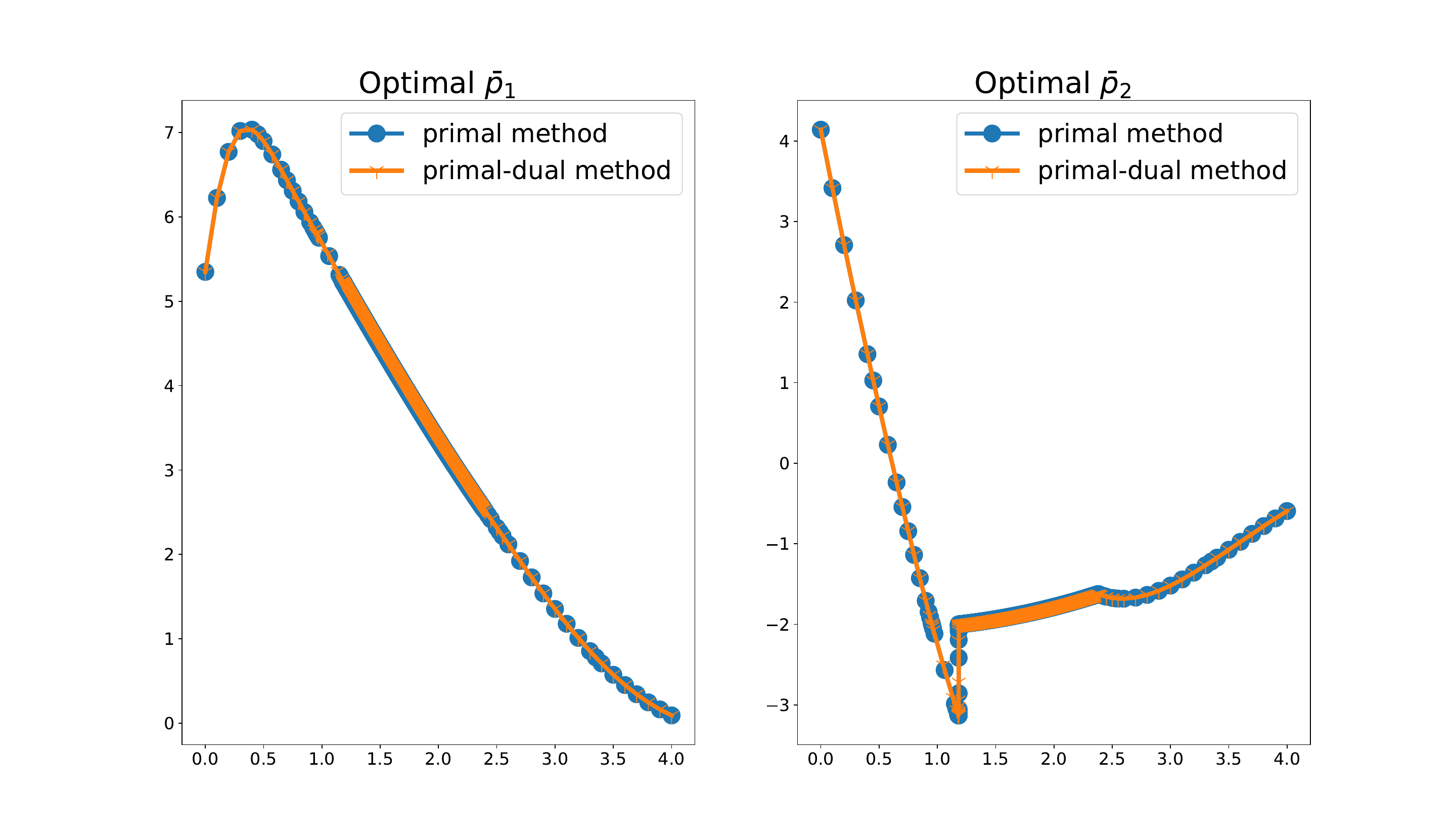}
    \caption{Optimal adjoint-states $p_1$ and $p_2$ using primal and primal-dual methods}
    \label{fig:vdp_adj_state}
\end{figure}

\begin{figure}[H]
    \centering
    \includegraphics[width=\textwidth]{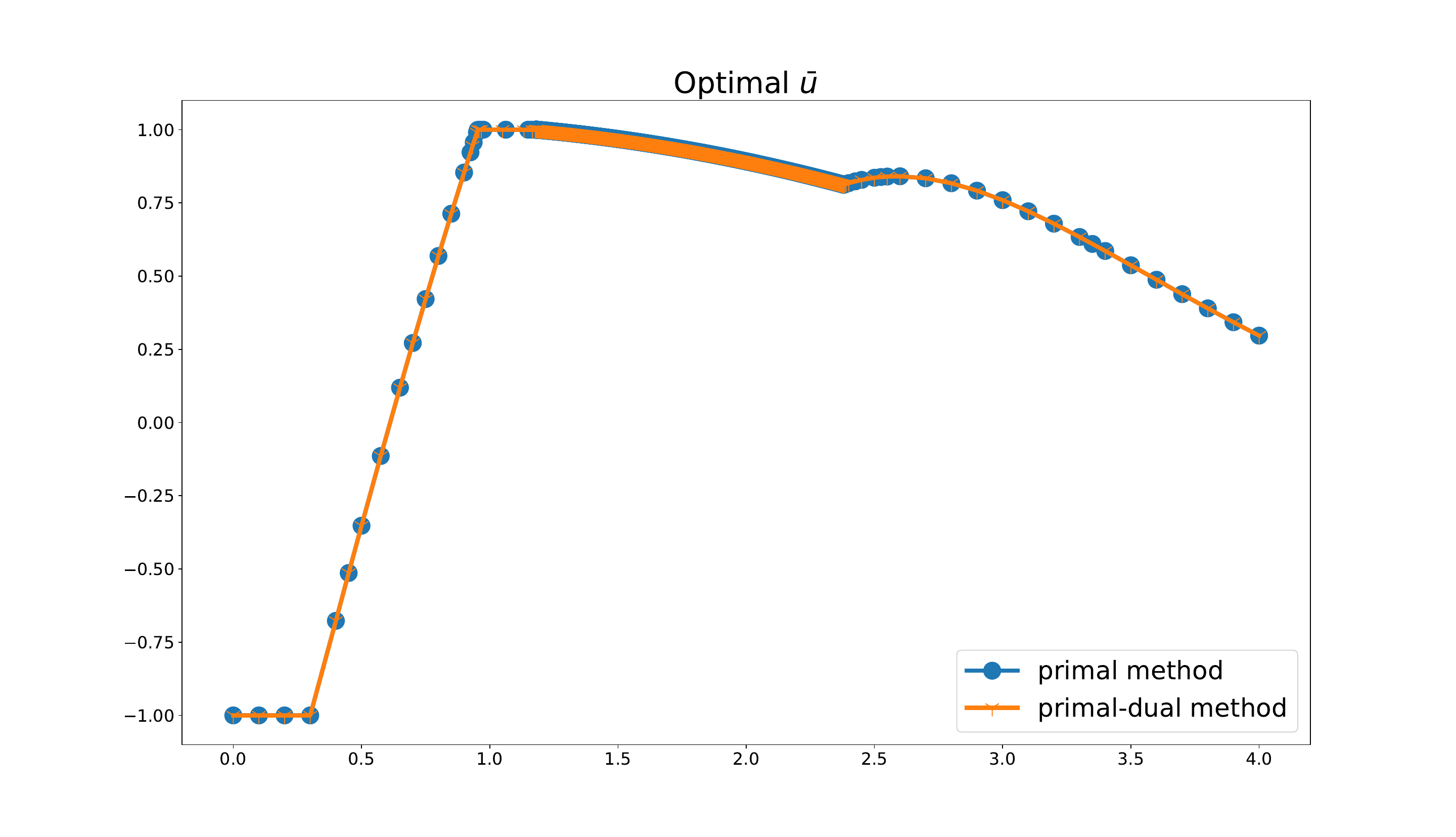}
    \caption{Optimal control $u$ using primal and primal-dual methods}
    \label{fig:vdp_control}
\end{figure}

First, \cref{fig:vdp_state,fig:vdp_adj_state,fig:vdp_control} show that both implementations converge to the same stationary point. However, \cref{tab:vdp} shows that the minimal decay ratio ensuring convergence is dramatically lower using a primal-dual implementation than the primal implementation. Therefore, the primal method achieves convergence in 16 $\epsilon$-iterations versus only 2 for the primal-dual method. \Cref{fig:vdp_state,fig:vdp_adj_state,fig:vdp_control,tab:vdp} also show that both methods produce close to identical numerical results. However, as displayed in \cref{tab:vdp} the primal-dual implementation is faster than the primal one.
\begin{center}
    \begin{table}[h!]
    \caption{Performance comparison between \cref{alg:primal} and \cref{alg:primal_dual} for the Van der Pol problem}
        \begin{tabular}{r|c|c|c|c}
        Method & decay ratio $\alpha$ & number of iterations & final length of time array & exec. time\\
          \hline
        Primal & 0.35 & 17 & 812 & 2.55 s\\
        Primal-dual & $10^{-7}$ & 2 & 797 & 1.92 s
        \end{tabular}
    \label{tab:vdp}
    \end{table}
\end{center}

\subsection{Zermelo's problem}
\subsubsection{Problem presentation and two-point boundary value problems}
The second example is the Zermelo's navigation problem \cite{CAILLAU2023559} which writes
\begin{equation*}
    \min_{(u,x,T)} \left[J(x,u,T) := T \right]
\end{equation*}
under the following constraints
\begin{align*}
    \dot{x}_1& = u_2  \cos(u_1) + h(x_2)\\
    \dot{x}_2 &= u_2 \sin(u_1)\\
    g(x_1,x_2)&\leq 0\\
    u_1 & \in [0, 2\pi]\\
    u_2 & \in [0, 1]\\
    x_1(0), x_2(0) & = 0\\
    x_1(T) & = 20\\
    x_2(T) & = 1
\end{align*}
with
\begin{align*}
    h(x_2) & := 3 + \frac{x_2(1 - x_2)}{5}\\
    g(x_1,x_2) &:=- \frac{(x_1 - 10)^2}{4} - \frac{(x_2-0.4)^2}{10^{-2}} +4
\end{align*}
The dynamics of this problem does not satisfy the \textbf{AP}. To transform this free final-time problem into a fixed end-time problem, we use the following change of variable $s := t/T$ and consider $T$ as a state variable satisfying $\dot{T} = 0$. In the following we denote $x_3:=T$. The fixed end-time problem is then
\begin{equation*}
    \min_{(x,u)} \left[ J(x,u) := x_3(1)\right]
\end{equation*}
under the following constraints
\begin{align*}
    \dot{x}_1& = x_3\left(u_2  \cos(u_1) + h(x_2)\right)\\
    \dot{x}_2 &= x_3 u_2 \sin(u_1)\\
    \dot{x}_3 & = 0\\
    g(x_1,x_2)&\leq 0\\
    u_1 & \in [0, 2\pi]\\
    u_2 & \in [0, 1]\\
    x_1(0), x_2(0) & = 0\\
    x_1(1) & = 20\\
    x_2(1) & = 1
\end{align*}
 One can easily see that any primal solution ${\rm S_P}(\epsilon)$ is solution of the following TPBVP
 \begin{align*}
    \dot{x}_1& = x_3\left(u_2  \cos(u_1) + h(x_2)\right)\\
    \dot{x}_2 &= x_3 u_2 \sin(u_1)\\
    \dot{x}_3 & = 0\\
    \dot{p}_1 & = \epsilon\frac{g'_{x_1}(x_1,x_2)}{g(x_1,x_2)}\\
    \dot{p}_2 & = - p_1 x_3 h'(x_2) +  \epsilon\frac{g'_{x_2}(x_1,x_2)}{g(x_1,x_2)}\\
    \dot{p}_3 & =- p_1 (u_2 \cos(u_1) + h(x_2)) - p_2 u_2 \sin(u_1)\\
     0 & = - p_1 x_3 u_2 \sin(u_1) + p_2 x_3 u_2 \cos(u_1) + \frac{\epsilon}{2\pi - u_1} - \frac{\epsilon}{u_1}\\
    0 & = p_1 x_2 \cos(u_1) + p_2 x_3 \sin(u_1) + \frac{\epsilon}{1 - u_2} - \frac{\epsilon}{u_2}\\
    x_1(0), x_2(0) & = 0\\
    x_1(1) & = 20\\
    x_2(1) & = 1\\
    p_3(0), p_3(1) & = 0
\end{align*}

and that any primal-dual solution ${\rm S_{PD}}(\epsilon)$ is solution of 
 \begin{align*}
    \dot{x}_1& = x_3\left(u_2  \cos(u_1) + h(x_2)\right)\\
    \dot{x}_2 &= x_3 u_2 \sin(u_1)\\
    \dot{x}_3 & = 0\\
    \dot{p}_1 & = -\lambda^g g'_{x_1}(x_1, x_2)\\
    \dot{p}_2 & = - p_1 x_3 h'(x_2)-\lambda^g g'_{x_2}(x_1, x_2)\\
    \dot{p}_3 & =- p_1 (u_2 \cos(u_1) + h(x_2)) - p_2 u_2 \sin(u_1)\\
     0 & = - p_1 x_3 u_2 \sin(u_1) + p_2 x_3 u_2 \cos(u_1) +\lambda^c_1 - \lambda^c_2\\
    0 & = p_1 x_2 \cos(u_1) + p_2 x_3 \sin(u_1) +\lambda^c_3 - \lambda^c_4\\
    0 & = {\rm FB}(\lambda^g, g(x_1, x_2), \epsilon)\\
    0 & = {\rm FB}(\lambda^c_1, u_1 - 2\pi, \epsilon)\\
    0 & = {\rm FB}(\lambda^c_2, -u_1, \epsilon)\\
    0 & = {\rm FB}(\lambda^c_3, u_2 - 1, \epsilon)\\
    0 & = {\rm FB}(\lambda^c_4, -u_2, \epsilon)\\
    x_1(0), x_2(0) & = 0\\
    x_1(1) & = 20\\
    x_2(1) & = 1\\
    p_3(0), p_3(1) & = 0
\end{align*}

\subsubsection{Initialization}
As illustrated on the left-hand side of \cref{fig:zermelo_state}, the first and second states $x_1, x_2$ are initialized differently depending on the solving algorithm. One can see that using a primal-dual method alleviates the constraint of strictly interior trajectory initialization. This is a strong advantage of this method since finding an interior trajectory might be challenging. The other variables are initialized with the same values for both methods.
\begin{align*}
    {\rm time}&\gets \begin{bmatrix}0,& 0.01,&\dots,&0.99,&1 \end{bmatrix} \\
    x_3&\gets \begin{bmatrix}20,&\dots,&20 \end{bmatrix}  \\
    p_1, p_2&\gets \begin{bmatrix}0,&\dots,&0 \end{bmatrix}\\
    p_3&\gets \begin{bmatrix}1,&\dots,&1 \end{bmatrix}\\
    u_1 & \gets \begin{bmatrix}\pi/2,&\dots,&\pi/2 \end{bmatrix}\\
    u_2 & \gets \begin{bmatrix}0.5,&\dots,&0.5 \end{bmatrix}\\
    \epsilon_0 & \gets 0.1
\end{align*}
The constraints multipliers for the primal-dual algorithm are initialized as follows
\begin{equation*}
    \lambda^g,\lambda^c_1,\dots, \lambda^c_4 \gets \begin{bmatrix}0,&\dots,&0 \end{bmatrix}
\end{equation*}

\subsubsection{Results}
\begin{figure}[H]
    \centering
    \includegraphics[width=\textwidth]{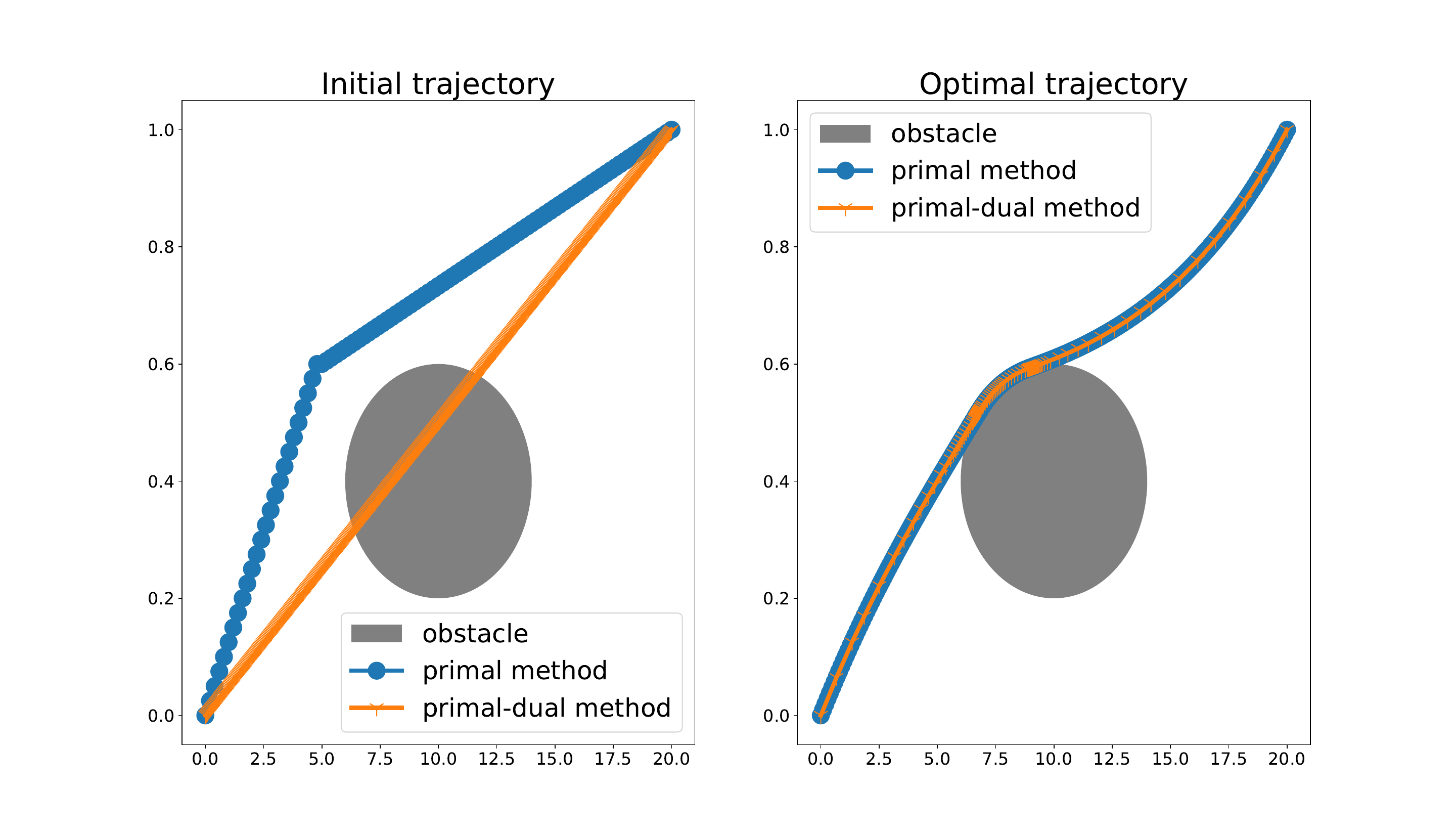}
    \caption{The figure on the left-hand side illustrates the initial trajectory $(x_1, x_2)$ used for both optimization methods. One can see that the initial primal-dual trajectory does not satisfy the state-constraint. The figure on the right-hand side illustrates the optimal trajectories from both methods.}
    \label{fig:zermelo_state}
\end{figure}

\begin{figure}[H]
    \centering
    \includegraphics[width=\textwidth]{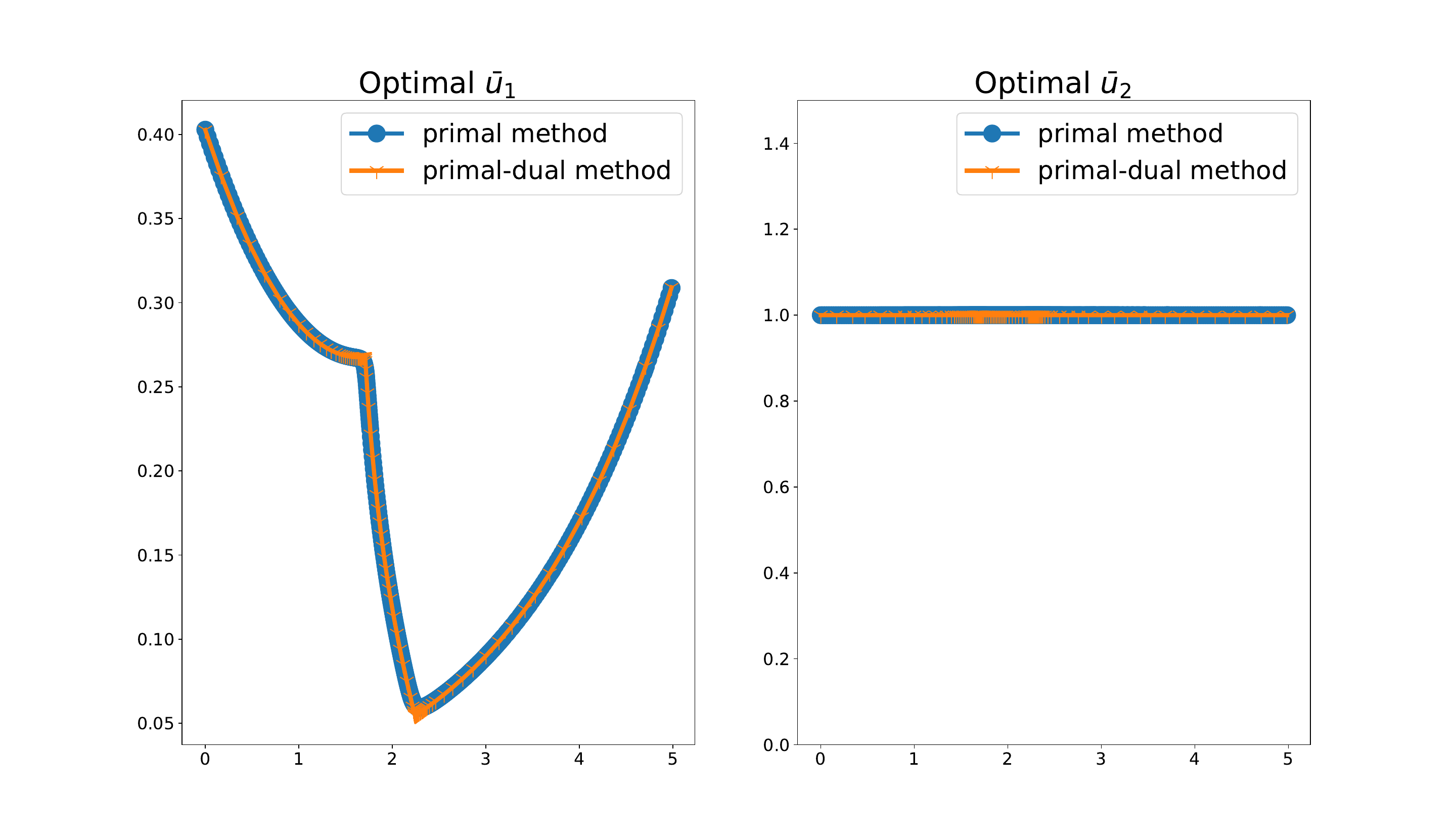}
    \caption{Optimal controls $u_1$ and $u_2$ using primal and primal-dual methods}
    \label{fig:zermelo_control}
\end{figure}

\begin{figure}[H]
    \centering
    \includegraphics[width=\textwidth]{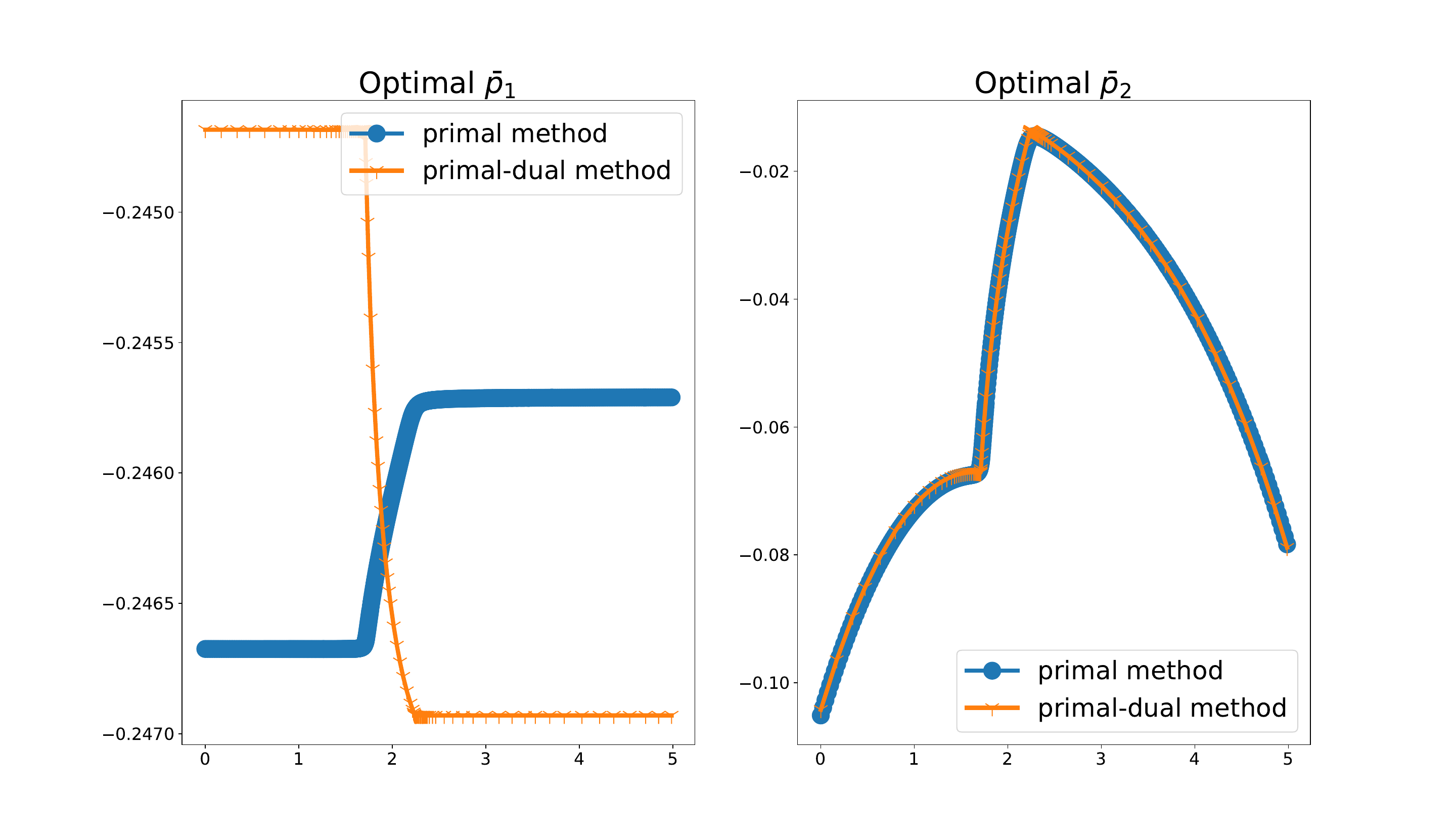}
    \caption{Optimal adjoint-states $p_1$ and $p_2$ using primal and primal-dual methods}
    \label{fig:zermelo_adj_state}
\end{figure}

\begin{figure}[H]
    \centering
    \includegraphics[width=\textwidth]{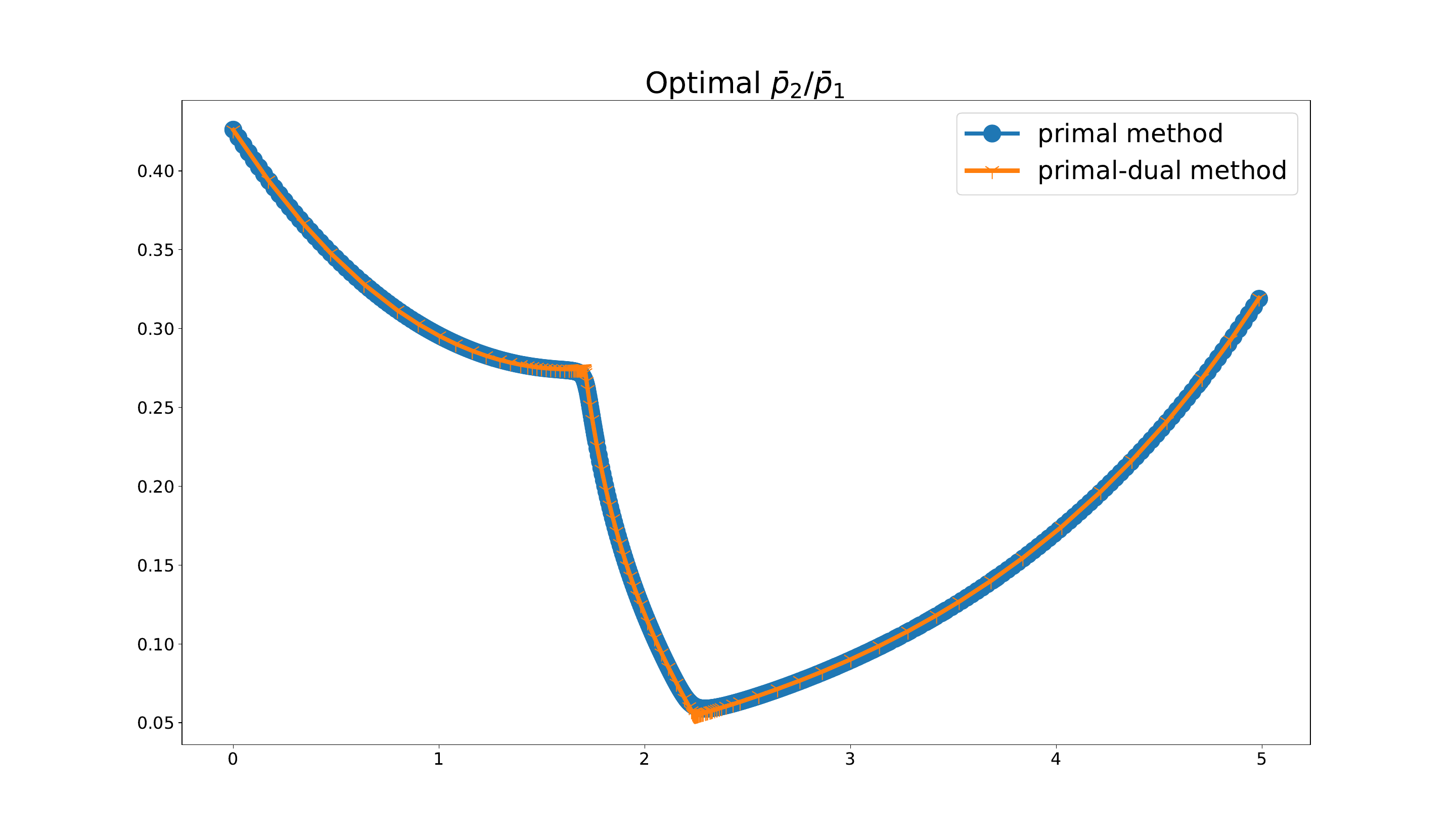}
    \caption{Optimal switching function $p_2 / p_1$ using primal and primal-dual methods}
    \label{fig:zermelo_switch_state}
\end{figure}

On the right-hand side of \cref{fig:zermelo_state}, one can see that both methods converge to the same optimal trajectory and avoid the obstacle by passing over it. As expected for a minimum time problem, one can see on the right-hand side of \cref{fig:zermelo_control} that the optimal speed $u_2$ is, as expected, constantly saturated at its maximal value. In addition, on the left-hand side of \cref{fig:zermelo_control}, one can see that the optimal navigation angle $u_1$ is the same for both methods, and it also appears that the mesh refinement produces a much finer grid for the primal implementation compared to the primal-dual one along the state-constrained arc. This finer grid yields more mesh points for the primal implementation than the primal-dual, as displayed in \cref{tab:zermelo}. The greater number of mesh points together with a greater decay ratio directly translates into a much longer execution time for the primal method compared to the primal-dual one. On the left-hand side of \cref{fig:zermelo_adj_state}, one can see that the first adjoint-state $p_1$ differs depending on the solving algorithm. However, as described in \cite{CAILLAU2023559}, the optimal control depends on the ratio $p_2 / p_1 $. One can see in \cref{fig:zermelo_switch_state} that the ratio $p_2 / p_1$ is the same for both methods; hence, the identical optimal controls as displayed in \cref{fig:zermelo_control}.
\begin{center}
    \begin{table}[h!]
        \caption{Performance comparison between \cref{alg:primal} and \cref{alg:primal_dual} for Zermelo's problem}
            \begin{tabular}{r|c|c|c|c}
            Method & decay ratio $\alpha$ & number of iterations & final length of time array & exec. time\\
              \hline
            Primal & 0.9 & 82 & 496 & 34.83 s\\
            Primal-dual & 0.5 & 21 & 132 & 4.99 s
            \end{tabular}
    \label{tab:zermelo}
    \end{table}
\end{center}

\subsection{Goddard's problem}
\subsubsection{Problem presentation and two-point boundary value problems}
The final example is the Goddard's problem \cite{seywald1993goddard} which writes
\begin{equation*}
    \min_{(x, u, T)}\left[J(x,u,T) :=\int_0^T -v \xd t\right]
\end{equation*}
under the following constraints
\begin{align*}
    \dot{h} & = v\\
    \dot{v} & = \frac{u - d(h,v)}{m} - \frac{1}{h^2}\\
    \dot{m} & =  - 2u\\
    h(0),\; m(0) & = 1\\
    v(0) & = 0\\
    m(T) &= 0.6\\
    q(h,v) & \leq 0\\
    u & \in [0, 3.5]
\end{align*}
with
\begin{align*}
    d(h,v) &:= 310 v^2\textrm{exp}(500(1-h))\\
    q(h,v) &:= 20 d(v,h) - 10
\end{align*}
As we did for the Zermelo's problem, we transform this free end-time problem in the following fixed end-time problem.
\begin{equation*}
    \min_{(x, u)}\left[J(x,u) :=\int_0^1 -T v \xd t\right]
\end{equation*}
under the following constraints
\begin{align*}
    \dot{h} & = Tv\\
    \dot{v} & = T\left(\frac{u - d(h,v)}{m} - \frac{1}{h^2}\right)\\
    \dot{m} & =  - 2Tu\\
    \dot{T} & = 0\\
    h(0),\; m(0) & = 1\\
    v(0) & = 0\\
    m(1) &= 0.6\\
    q(h,v) & \leq 0\\
    u & \in [0, 3.5]
\end{align*}
This problem is challenging because the optimal control trajectory consists of a sequence made of a bang-bang arc, a singular arc, a state-constrained arc, and finally, a bang-bang arc. One can easily verify that any primal solution ${\rm S_P}(\epsilon)$ is solution of the following TPBVP
\begin{align}
    \dot{h} & = Tv\nonumber\\
    \dot{v} & = T\left(\frac{u - d(h,v)}{m} - \frac{1}{h^2}\right) \nonumber\\
    \dot{m} & =  - 2Tu\nonumber\\
    \dot{T} & = 0\nonumber\\
    \dot{p}_h & = -\left(\frac{2}{h^3}-\frac{d'_h(h,v)}{m}\right)Tp_v+\epsilon\frac{q'_h(h,v)}{q(h,v)}\nonumber\\
    \dot{p}_v & = T\left(1-p_h+\frac{d'_v(h,v)}{m}p_v\right)+\epsilon\frac{q'_v(h,v)}{q(h,v)} \nonumber\\
    \dot{p}_m & = \frac{u-d(h,v)}{m^2}Tp_v \nonumber\\
    \dot{p}_T &= v-p_hv- \left(\frac{u- d(h,v)}{m}-\frac{1}{h^2} \right)p_v +2up_m \nonumber\\
    0 &=\frac{Tp_v}{m}-2Tp_m +\epsilon\left(1/(3.5-u) -1/u\right)\label{eq:goddard_Huprime_SP}\\
    h(0),\; m(0) & = 1\nonumber\\
    v(0),\; p_T(0),\; p_h(1),\; p_v(1),\; p_T(1) & = 0\nonumber\\
    m(T) &= 0.6\nonumber
\end{align}

and that any primal-dual solution ${\rm S_{PD}}(\epsilon)$ is solution of 
\begin{align}
    \dot{h} & = Tv\nonumber\\
    \dot{v} & = T\left(\frac{u - d(h,v)}{m} - \frac{1}{h^2}\right) \nonumber\\
    \dot{m} & =  - 2Tu\nonumber\\
    \dot{T} & = 0\nonumber\\
    \dot{p}_h & = -\left(\frac{2}{h^3}-\frac{d'_h(h,v)}{m}\right)Tp_v-\lambda^g q'_h(h,v)\nonumber\\
    \dot{p}_v & =T\left(1-p_h+\frac{d'_v(h,v)}{m}p_v\right)-\lambda^g q'_v(h,v) \nonumber\\
    \dot{p}_m & = \frac{u-d(h,v)}{m^2}Tp_v \nonumber\\
    \dot{p}_T &= v-p_hv- \left(\frac{u- d(h,v)}{m}-\frac{1}{h^2} \right)p_v +2up_m \nonumber\\
    0 &=\frac{Tp_v}{m}-2Tp_m +\lambda^c_1 - \lambda^c_2\label{eq:goddard_Huprime_SPD}\\
    0 & = {\rm FB}(\lambda^g, q(h,v),\epsilon)\nonumber\\
    0 & = {\rm FB}(\lambda^c_1, u - 3.5,\epsilon) \nonumber\\
    0 & = {\rm FB}(\lambda^c_2, -u,\epsilon) \nonumber \\
    h(0),\; m(0) & = 1\nonumber\\
    v(0),\; p_T(0),\; p_h(1),\; p_v(1),\; p_T(1) & = 0\nonumber\\
    m(T) &= 0.6\nonumber
\end{align}
\subsubsection{Initialization}
For both algorithms the initial values of time, state, adjoint-states, control and the penalty parameter are as follows

\begin{align*}
    {\rm time}&\gets \begin{bmatrix}0, &0.01, &\dots, &0.99, &1\end{bmatrix}\\
    h & \gets \begin{bmatrix}1.2, &\dots,  &1.2\end{bmatrix}\\
    v & \gets \begin{bmatrix}0.05, &\dots, &0.05\end{bmatrix}\\
    m & \gets \begin{bmatrix}1, &\dots, &1\end{bmatrix}\\
    T & \gets \begin{bmatrix}0.3, &\dots, &0.3\end{bmatrix}\\
    p_h, p_m, p_T & \gets \begin{bmatrix}0, &\dots, &0\end{bmatrix}\\
    p_v & \gets \begin{bmatrix}1, &\dots, &1\end{bmatrix}\\
    \epsilon_0 & \gets 0.1
\end{align*}
And the constraints multipliers for the primal-dual algorithm are initialized as follows
\begin{equation*}
    \lambda^g,\lambda^c_1,\lambda^c_2 \gets \begin{bmatrix}0, &\dots, &0\end{bmatrix}
\end{equation*}

\subsubsection{Results}
\begin{figure}[H]
    \centering
    \includegraphics[width=\textwidth]{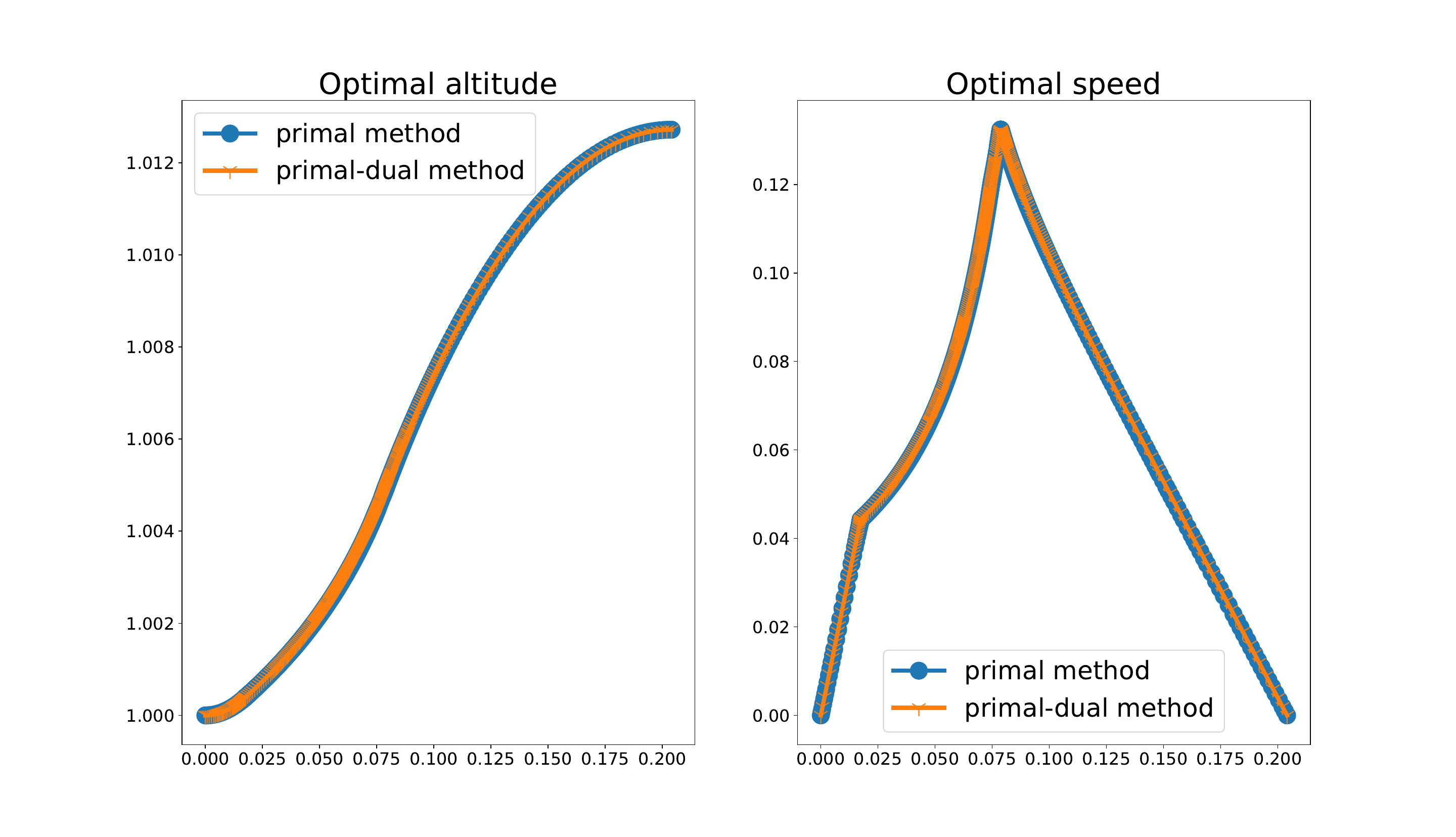}
    \caption{Optimal altitude $h$ and speed $v$ for primal and primal dual methods}
    \label{fig:goddard_hv}
\end{figure}

\begin{figure}[H]
    \centering
    \includegraphics[width=\textwidth]{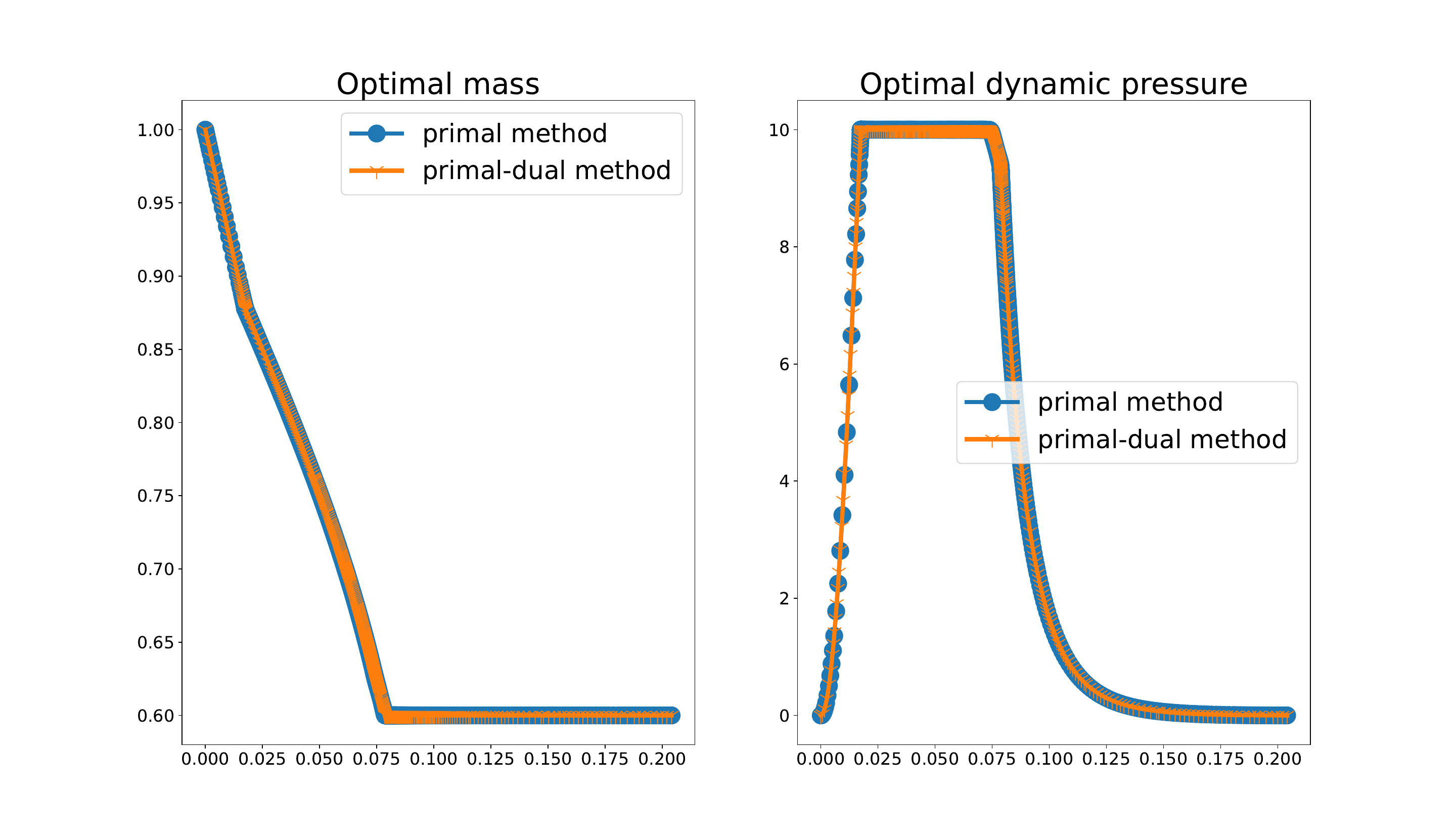}
    \caption{Optimal mass $m$ and dynamic pressure $q$ for primal and primal dual methods}
    \label{fig:goddard_mq}
\end{figure}

\begin{figure}[H]
    \centering
    \includegraphics[width=\textwidth]{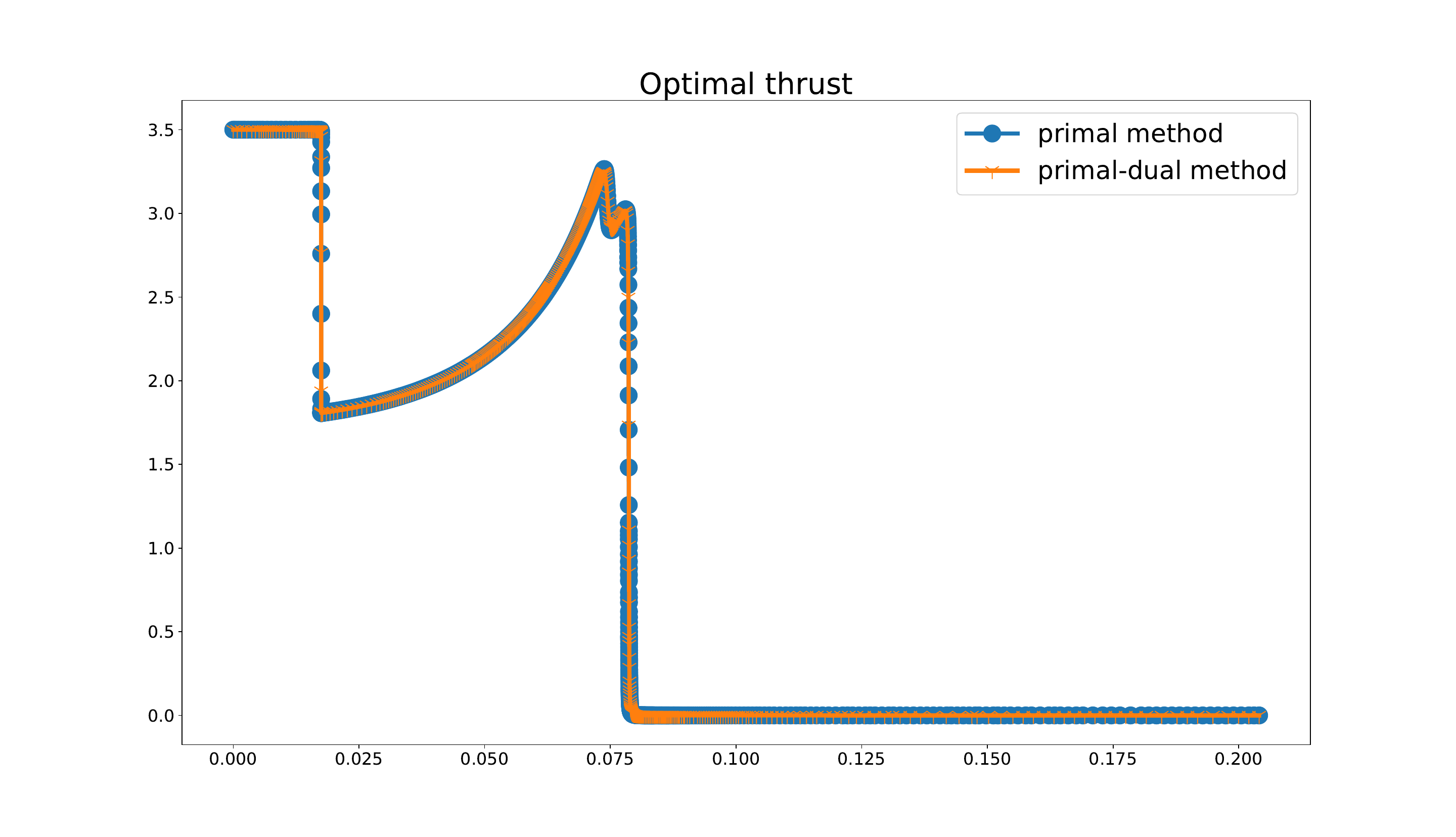}
    \caption{Optimal thrust $u$ for primal and primal dual method}
    \label{fig:goddard_u}
\end{figure}

The optimal altitude and speed corresponding to each solving algorithm are displayed in \cref{fig:goddard_hv}, the optimal mass and dynamic pressure are displayed in \cref{fig:goddard_mq}, and the optimal thrust is displayed in \cref{fig:goddard_u}. Both solving algorithms converge to the same optimal solution. Both algorithms allow for finding the well-known succession of bang-bang, singular, and state-constrained arcs for the optimal thrust. Finally, in \cref{tab:goddard}, one can see that the primal-dual method achieves better performance compared to the primal one in terms of decay ratio, and execution time.

\begin{center}
    \begin{table}[h!]
    
    \caption{Performance comparison between \cref{alg:primal} and \cref{alg:primal_dual} for Goddard's problem}
        \begin{tabular}{r|c|c|c|c}
        Method & decay ratio $\alpha$ & number of iterations & final length of time array & exec. time\\
          \hline
        Primal & 0.6 & 29 & 722 & 16.14 s\\
        Primal-dual & 0.25 & 11 & 501 & 4.11 s
        \end{tabular}
    \label{tab:goddard}
    \end{table}
\end{center}

\section{Data availability statement}
The numerical examples and the Differential Algebraic Equations (DAEs) solver used in \ref{sec:algoExample} is available online in a Gitlab repository: \url{https://ifpen-gitlab.appcollaboratif.fr/detocs/ipm_ocp}. The solver is a two point boundary differential algebraic equations solver adapted from \cite{Kierzenka2001ABS} to solve index-1 differential algebraic equations.

\appendix

\section{Proof of \cref{thm:interior_all}}
\label{sec:proof_lemma_siam}
\subsection{Technical results and definition.}
To prove \cref{thm:interior_all}, we need some preliminary results stated below in \cref{prop:lower_bound_measure}, \cref{prop:maximal_distance}, and the following definition.
\begin{dfntn}[State-constraint measure]
\label{def:mui}
For all $(u, x^0) \in \U \times \R^n$, and for all $E \subset \R $, we denote $g_i(x[u,x^0])\#\meas$ the push-forward $g_i$-measure of $E$ defined as follows
\begin{equation}
    \label{eq:def_mui_alpha1_alpha2}
   g_i(x[u,x^0])\#\meas(E) := \meas\left(\left(g_i\circ x[u, x^0]\right)^{-1}(E)\right)
\end{equation}
\end{dfntn}

\begin{prpstn}
\label{prop:lower_bound_measure}
Let $K<+\infty$, for all $(u, x^0) \in B_{\norm{.}_{\Vad}}(0,K)$, and for all Lebesgue-measurable set $E \subseteq g_i\circ x[u, x^0]([0,T])$, the state-constraint measure from \cref{def:mui} is lower bounded as follows
\begin{equation}
\label{eq:lower_bound_measure}
   g_i(x[u,x^0])\#\meas(E)\geq \const(f,g)\meas(E)
\end{equation}
\end{prpstn}

\begin{proof}
Assume that \cref{ass:bounded_control_set} holds, then for all $(u, x^0) \in B_{\norm{.}_{\Vad}}(0,K)$ we have
\begin{align*}
    \vert g_i(x[u, x^0](t)) - g_i(x[u, x^0](s))\vert &\leq \const(g)\norm{x[u, x^0](t) - x[u, x^0](s)}\\
    &= \const(g)\norm{ \int_s^t f(x[u, x^0](\tau), u(\tau))}\xdt\\
    &\leq \const(f,g)\vert t-s\vert
\end{align*}
To prove the proposition, we only need to prove the lower bound holds on any interval $(\alpha_1,\alpha_2)\subseteq E$. From the continuity of $g_i$, $\exists t_1,t_2$ such that $g_i(x[u, x^0](t_1))=\alpha_1$, $g_i(x[u, x^0](t_2))=\alpha_2$ and such that $(t_1,t_2)\subseteq g_i(x[u, x^0])^{-1}((\alpha_1,\alpha_2))$ and
\begin{equation*}
  g_i(x[u,x^0])\#\meas((\alpha_1,\alpha_2)) \geq \vert t_1 - t_2\vert \geq\const(f,g) \vert g_i(x[u, x^0](t_1)) - g_i(x[u, x^0](t_2))\vert \geq \const(f,g)\vert \alpha_1 - \alpha_2\vert  
\end{equation*}

\end{proof}

\begin{prpstn}
\label{prop:maximal_distance}
For all $\delta>0$, there exists $G_\delta, C_\delta\in (0,+\infty)$ such that for all $(u, x^0) \in \Vad_\infty$, there exists $(v,y^0) \in B_{\norm{.}_{\Vad}}((u, x^0),\delta)\cap\Vad$ satisfying the following conditions
\begin{align*}
    \sup_t g_i(x[v, y^0](t))&\leq -2G_\delta,\;i=1,\dots,n_g\\
    \esssup_t c_i(x[v, y^0](t), v(t))& \leq -2C_\delta,\;i=1,\dots,n_c
\end{align*}
and we also have
\begin{align*}
    g_i(x[v, y^0](t))\leq g_i(x[u, x^0](t))-G_\delta,\;\forall t\in S^g_{u,x^0}(G_\delta),\;i=1,\dots,n_g\\ 
    c_i(x[v, y^0](t), v(t))\leq c_i(x[u, x^0](t), u(t))-C_\delta,\;\textrm{a.e. } t\in S^c_{u,x^0}(C_\delta),\;i=1,\dots,n_c
\end{align*}
\end{prpstn}
\begin{proof}
Let $\delta>0$, and for all $(u, x^0)\in \Vad_\infty$ let us denote
\begin{equation*}
    \gamma_\delta(u, x^0) := \inf_{v\in B_{\norm{.}_{\Vad}}((u,x^0),\delta)\cap \Vad}\left\lbrace \sup_t g(x[v, y^0](t))\right\rbrace
\end{equation*}
From \cref{ass:interior_accessibility}, we have $\gamma_\delta(u, x^0)<0$. Then, $\forall (u, x^0)\in \Vad_\infty$, $\exists (v,y^0)\in B_{\norm{.}_{\Vad}}((u,x^0),\delta)\cap \Vad$ such that
\begin{equation*}
\sup_t g(x[v, y^0](t))\leq \gamma_\delta(u, x^0)\leq \sup_{(u, x^0)\in \Vad_\infty}\gamma_\delta (u, x^0):=-2G_\delta <0
\end{equation*}
In addition,  if $\supp^g_{u,x^0}(G_\delta)\neq \emptyset$ , then $\forall t\in \supp^g_{u,x^0}(G_\delta)$ we have
\begin{align*}
    g(x[v, y^0](t)) - g(x[u, x^0](t))\leq -2G_\delta + G_\delta =- G_\delta
\end{align*}
Now, let us denote  
\begin{equation*}
    \kappa_\delta(u,x^0) := \inf_{v\in B_{\norm{.}_{\Vad}}((u,x^0),\delta)\cap \Vad}\left\lbrace \esssup_t c(x[v, y^0](t), v(t))\right\rbrace
\end{equation*}
From \cref{ass:interior_accessibility}, for all $(u, x^0)\in \Vad_\infty$, $\exists (v_n, y^0_n)_n\in \VadStrict(n)$ converging to $(u, x^0)$. Thus, for all $\delta >0$, $\exists N_\delta[u, x^0]>0$, such that $\forall n\geq N_\delta[u, x^0]$, we have $(v_n,y_n^0)\in B_{\norm{.}_{\Vad}}((u, x^0), \delta)$, which yields $\kappa_\delta(u,x^0) \leq -1/N_\delta(u, x^0)$. Now, $\forall (u, x^0)\in \Vad_\infty$, $\exists (v,y^0)\in B_{\norm{.}_{\Vad}}((u,x^0),\delta)\cap \Vad$ such that
\begin{equation*}
    \esssup_t c(x[v, y^0](t), v(t))\leq \kappa_\delta(u, x^0)\leq \sup_{(u, x^0)\in \Vad_\infty}\kappa_\delta (u, x^0)
    =\sup_{(u, x^0)\in \Vad_\infty} -\frac{1}{N_\delta[u, x^0]}:=-2C_\delta <0
\end{equation*}
In addition,  if $\supp^c_{u,x^0}(C_\delta)\neq \emptyset$ , then $\forall t\in \supp^c_{u,x^0}(C_\delta)$ we have
\begin{equation*}
    c(x[v, y^0](t), v(t)) - c(x[u, x^0](t), u(t))\leq -2C_\delta + C_\delta =- C_\delta
\end{equation*}
\end{proof}

\label{sec:proof_thm_int_all}
\subsection{Proof of \cref{eq:interiority_state,eq:boundedness_lambda_g}}
\label{sec:proof_g_lg_boundedness}
To prove the result, it is sufficient to prove for $n_g = 1$, that is to say, for just one state constraint.Now, let us consider the following state-penalized OCP
\begin{subequations}
\label{eq:all_state_penalized_problem}
    \begin{align}
          \min_{(x, u)} J^1_\epsilon(x, u) &:=\varphi(x(T)) + \int_0^T\left[\ell(x(t),u(t))+\epsilon\sum_{i=1}^{n_g} \psi \circ g_i(x(t))\right]\xdt\label{eq:def_state_pen_problem_1}\\
          \dot{x}(t) & = f(x(t),u(t))\label{eq:def_state_pen_problem_2}\\
          h(x(0), x(T)) & = 0 \label{eq:def_state_pen_problem_3}\\
          c(x(t), u(t))& \leq 0 \label{eq:def_state_pen_problem_4}
    \end{align}
\end{subequations}
\begin{lmm}
\label{thm:interior_state}
For all $\epsilon \in(0,\epsilon_0)$, any associated locally optimal solution $(x[\ue, x^0_\epsilon], \ue)$ of Problem \cref{eq:all_state_penalized_problem} is such that \cref{eq:interiority_state,eq:boundedness_lambda_g} hold
\end{lmm}
\begin{proof}
Assume $(\ue, x^0_\epsilon) \in \Vad_\infty$ is a locally optimal solution of Problem \cref{eq:all_state_penalized_problem} satisfying 
\begin{equation}
\label{eq:ass_touchpoint_g}
    \sup_t g(x[\ue, x_\epsilon^0](t)) = 0
\end{equation}
From \cref{prop:maximal_distance}, $\forall \delta >0, \exists (\vd, \xdelta^0)\in B_{\norm{.}_{\Vad}}((\ue, x^0_\epsilon),\delta) \cap \Vad$ and $G_\delta >0$ such that
\begin{subequations}
    \begin{align*}
        g(x[\vd, \xdelta^0](t))&\leq -2 G_\delta,\; \forall t \in[0,T]\\
         g(x[\vd, \xdelta^0](t))&\leq g(x[\ue, x^0_\epsilon](t))-G_\delta,\; \forall t\in \supp^g_{\ue, x^0_\epsilon}(G_\delta) \label{eq:inequality_sur_suppbue}
    \end{align*}
\end{subequations}
with $\supp^g_{\ue, x^0_\epsilon}(G_\delta)\neq \emptyset$. In the following, to alleviate the notations, we denote 
\begin{align*}
    \zd &:= (\vd, \xdelta^0)\\
    \ze &:= (\ue, x^0_\epsilon)\\
    \Delta z &:= \zd - \ze\\
    \Delta x &:= x[\zd] - x[\ze]\\
    \Delta g &:= g(x[\zd]) - g(x[\ze])
\end{align*}
Now, one can exhibit an upper-bound on the difference $J^1_\epsilon(x[\zd], \vd)-J^1_\epsilon(x[\ze], \ue)$ as follows
\begin{equation*}
    J^1_\epsilon(x[\zd], \vd)-J^1_\epsilon(x[\ze], \ue) := \Delta_1 + \epsilon \Delta_2\label{eq:def_split_g_epsilon}
\end{equation*}
where
\begin{align*}
    \Delta_1 &:= 
        \varphi(x[\zd](T))-\varphi(x[\ze](T)) 
        + \int_0^T\left[\ell(x[\zd](t),\vd(t))-\ell(x[\ze](t),\ue(t))\right]\xdt\\
     \Delta_2 &:= \int_0^T\left[\psi\circ g(x[\zd](t))-\psi \circ g(x[\ze](t))\right]\xdt
\end{align*}
 Now, let us upper-bound $\Delta_1$. From \cref{ass:bounded_if_state_set}, $x_\epsilon(0)$ is bounded. From \cref{ass:bounded_control_set,eq:uniform_bounded_u}, $u_\epsilon$ is  $\xLinfty$-bounded, thus $\xLone$-bounded and there exists $R_v<+\infty$, such that $z_\epsilon, z_\delta \in B_{\norm{.}_{\Vad}}(0,R_v)$. Using \cref{ass:bounded_control_set,eq:ass_bibo}, we have
\begin{align}
    \Delta_1 \leq&
    \int_0^T\const(\ell)\left(\norm{x[\zd](t)-x[\ze](t)}+\norm{\vd(t)-\ue(t)}\right)\xdt
    +\const(\varphi)\norm{\Delta x}_{\xLinfty}\nonumber\\
    \leq&\const(\ell, T, \varphi)\norm{\Delta x}_{\xLinfty}+\const(\ell)\norm{\vd-\ue}_{\xLone} \leq\const(\ell,f,g, \varphi, T, R_v, R_x)\label{eq:upperbound_delta1}
\end{align}
Now, let us upper-bound $\Delta_2$. To do so, let us introduce the following useful subsets of $[0,T]$
 \begin{align}
    E_1&:=(g\circ x[\ze])^{-1}\left((-\infty,-G_\delta]\right)\label{eq:def_set_a1}\\
    E_2(\rho)&:=(g\circ x[\ze])^{-1}\left((-G_\delta,-\rho]\right)\label{eq:def_set_a2}\\
    E_3(\rho)&:=(g\circ x[\ze])^{-1}\left((-G_\delta,-\rho)\right)\label{eq:def_set_a3}
\end{align}
Given \cref{eq:ass_touchpoint_g}, for all $\rho\in [0,G_\delta)$, these sets are not empty and $\forall t \notin \left(E_1 \cup E_2(\rho)\right)$ we have $\psi\circ g(x[\zd](t))- \psi\circ g(x[\ze](t))<0$ which yields
\begin{equation*}
    \Delta_2 \leq \int_{E_1} \psi \circ g(x[\zd](t)) - \psi \circ g(x[\ze](t))\xdt
    + \int_{E_2(\rho)} \psi\circ g(x[\zd](t))- \psi \circ g(x[\ze](t))\xdt
\end{equation*}
By convexity of the $\log$-penalty, i.e. $\psi$, we have
\begin{equation*}
    \label{eq:majoS1}
    \int_{E_1} \psi\circ g(x[\zd](t)) - \psi\circ g(x[\ze](t))\xdt\leq
    \int_{E_1} \psi'(G_\delta) \norm{\Delta g}_{\xLinfty} \xdt:=\const(T,f,g,G_\delta)
\end{equation*}
In addition, we also have
\begin{equation*}
     \int_{E_2(\rho)}\psi\circ g(x[\zd](t)) - \psi\circ g(x[\ze](t))\xdt = \int_{E_2(\rho)} \left(\int_0^1 \psi'(g(x[\ze](t)) + s \Delta g(t)) \xd s\right)\Delta g(t) \xdt
\end{equation*}
Since $\forall t\in E_2(\rho), \,\Delta g(t)  <  -G_\delta$, we also have 
\begin{equation}
    \label{eq:majoS2}
    \int_{E_2(\rho)}\psi\circ g(x[\zd](t)) - \psi\circ g(x[\ze](t))\xdt \leq -G_\delta\int_{E_2(\rho)}\left(\int_{0}^{1} \psi'(g(x[\ze](t)) + s\Delta g(t))\xd s\right) \xdt
\end{equation}
From the mean value theorem, $\forall t\in E_2(\rho)$, $\exists \sigma_t\in(0,1)$ such that
\begin{equation}
\label{eq:mean_value_g}
    \psi'(g(x[\ze](t)) + \sigma_t \Delta g(t))= \int_{0}^{1} \psi'(g(x[\ze](t)) + s\Delta g(t))\xd s
\end{equation}

and we have
\begin{equation*}
    \psi'\circ g(x[\zd](t)) < \psi'(g(x[\ze](t)) + \sigma_t \Delta g(t)) <\psi'\circ g(x[\ze](t))
\end{equation*}
From the intermediate value theorem, $\exists \bar{\sigma}\in(0,1)$ such that $\forall t\in E_2(\rho)$ we have
\begin{equation}
\label{eq:lower_bound_intermediate}
    \psi'(g(x[\ze](t)) + \sigma_t \Delta g(t)) \geq (1 - \bar{\sigma})\psi'\circ g(x[\ze](t)) + \bar{\sigma} \psi'\circ g(x[\zd](t))
\end{equation}
Gathering \cref{eq:majoS2,eq:mean_value_g,eq:lower_bound_intermediate} yields
\begin{multline}
    \int_{E_2(\rho)}\psi\circ g(x[\zd](t)) - \psi\circ g(x[\ze](t))\xdt \leq - G_\delta \int_{E_2(\rho)}\bigg( (1-\bar{\sigma})\psi' \circ g(x[\ze](t)) + \bar{\sigma} \psi' \circ g(x[\zd](t)) \bigg)\xdt\\
    \leq- G_\delta\left( (1 - \bar{\sigma})\int_{E_2(\rho)}\psi' \circ g(x[\ze](t))\xdt + \const(\bar{\sigma}, \psi, T,g, G_\delta)\right)\label{eq:upperbound_delta2_over_E2}
\end{multline}
Gathering \cref{eq:majoS1,eq:upperbound_delta2_over_E2} we have 
\begin{equation}
    \label{eq:upperbound_delta2}
    \epsilon \Delta_2 \leq \const(T,g,G_\delta,\epsilon_0, \psi, \bar{\sigma}) -\epsilon G_\delta(1 - \bar{\sigma})\int_{E_2(\rho)}\psi' \circ g(x[\ze](t))\xdt
\end{equation}

Now, let us prove that any optimal solution is strictly interior with respect to the state constraint. The proof is by contradiction. Using \cref{def:mui}, one can make the following change in measure
\begin{equation*}
\int_{E_3(\rho)}\psi' \circ g(x[\ze](t))\xdt = \int_{-G_\delta}^{-\rho}\psi'(s)\xd g(x[\ze])\#\meas(s)
\end{equation*}
Then, using \cref{prop:lower_bound_measure,eq:def_set_a2,eq:def_set_a3} yields
\begin{align}
    \int_{E_2(\rho)}\psi' \circ g(x[\ze](t))\xdt&\geq\int_{E_3(\rho)}\psi' \circ g(x[\ze](t))\xdt \geq \const(f,g)\int_{-G_\delta}^{-\rho} \psi'(s)\xd s = \const(f,g)\left(\psi(-\rho)-\psi(-G_\delta)\right)\label{eq:proof_change_variable}
\end{align}
Gathering \cref{eq:def_split_g_epsilon,eq:upperbound_delta1,eq:upperbound_delta2,eq:proof_change_variable} yields that, for all $\rho>0$, we have
\begin{equation*}
    J^1_\epsilon(\zd) - J^1_\epsilon(\ze) \leq\const(\ell,f,g, \varphi, T, \epsilon_0,\psi, G_\delta, R_v,R_x, \bar{\sigma}) -\epsilon G_\delta \const(f,g, \bar{\sigma}) \left(\psi(-\rho)-\psi(-G_\delta)\right)
\end{equation*}
For $\rho$ small enough, this yields $J^1_\epsilon(\zd)<J^1_\epsilon(\ze)$ and contradicts the local optimality of $\ze$ and proves \cref{eq:interiority_state}. Thus, we have $(g(x[\ze]))^{-1}(\{0\}) =\emptyset$ which yields
\begin{equation}
    [0,T]\ = \lim_{\rho \rightarrow 0}(g(x[\ze]))^{-1}((-\infty, \rho))
    \label{eq:density_S3_S2}
\end{equation}
Hence, using \cref{def:mui}, one has
\begin{equation*}
    \norm{\psi'\circ g(x[\ze])}_{\xLone} =\int_0^T\psi' \circ g(x[\ze](t))\xdt:=\lim_{\rho\rightarrow 0}\int_{-\infty}^{-\rho}\psi'(s)\xd g(x[\ze])\#\meas(s)
\end{equation*}
Now, let us prove \cref{eq:boundedness_lambda_g} by contradiction and assume that
\begin{equation*}
    \forall K_g>0, \exists \epsilon>0 \st \norm{\epsilon \psi'\circ g(x[\ze]) }_{\xLone}>K_g
\end{equation*}
Then, from \cref{eq:def_set_a1,eq:def_set_a2,eq:density_S3_S2}, one has
\begin{equation}
\lim_{\rho \rightarrow 0}\epsilon\int_{-G_\delta}^{-\rho}\psi'(s)\xd g(x[\ze])\#\meas(s) > K_g - \epsilon \int_{-\infty}^{-G_\delta} \psi'(s) \xd g(x[\ze])\#\meas(s)> K_g -\frac{\epsilon_0 T}{G_\delta} \label{eq:upperbound_l1_gammaxprime}
\end{equation}
Gathering \cref{eq:upperbound_delta1,eq:upperbound_delta2,eq:upperbound_l1_gammaxprime} yields
\begin{equation*}
    \Delta_1 + \epsilon \Delta_2 \leq 
    \const(\ell,f,g, \varphi,T, \epsilon_0,\psi, G_\delta, R_v, R_x, \bar{\sigma})
-G_\delta(1 - \bar{\sigma}) K_g
\end{equation*}
Since $G_\delta(1 - \bar{\sigma})>0$, $\exists K_g>0$ such that $\Delta_1 + \epsilon \Delta_2<0$ which contradicts the optimality of $\ze$, proves \cref{eq:boundedness_lambda_g} and concludes the proof.
\end{proof}

\subsection{Proof of \cref{eq:boundedness_lambda_c}}
\begin{lmm}
\label{thm:interior_mixed_const}
There exists a constant $K_c<+\infty$ such that for all $\epsilon \in  (0,\epsilon_0)$ and for any $(x[\ue, x^0_\epsilon], \ue)$ locally optimal solution of Problem \cref{eq:def_log_barrier_ocp}, equation \cref{eq:boundedness_lambda_c} holds.
\end{lmm}

\begin{proof}
It is sufficient to prove the case where $n_c=1$, i.e., when there is a single mixed constraint. From \cref{prop:maximal_distance}, $\forall \delta >0, \exists (\vd, \xdelta^0)\in B_{\norm{.}_{\Vad}}((\ue, x^0_\epsilon),\delta) \cap \Vad$ and $C_\delta >0$ such that
\begin{subequations}
    \begin{align*}
        c(x[\vd, \xdelta^0](t), \vd(t))&\leq -2 C_\delta,\; \textrm{ a.a. }t\in[0,T]\\
         c(x[\vd, \xdelta^0](t), \vd(t))&\leq c(x[\ue, x^0_\epsilon](t), \ue(t))-C_\delta,\; \forall t\in \supp^c_{\ue, x^0_\epsilon}(C_\delta) 
    \end{align*}
\end{subequations}
with $\supp^c_{\ue, x^0_\epsilon}(C_\delta)\neq \emptyset$. In the following, to alleviate the notations, we denote 
\begin{align*}
    \zd &:= (\vd, \xdelta^0)\\
    \ze &:= (\ue, x^0_\epsilon)\\
    \Delta z &:= \zd - \ze\\
    \Delta x &:= x[\zd] - x[\ze]\\
    \Delta g &:= g(x[\zd]) - g(x[\ze])\\
    \Delta c &: = c(x[\zd],\vd) - c(x[\ze], \ue)
\end{align*}
In addition, From \cref{thm:interior_state}, and by continuity of the mapping $z\mapsto x[z]$ one can chose $\delta>0$ such that the following holds
\begin{align}
    \sup_t g(x[\zd](t))& < 0\label{eq:interiority_state_mixed}\\
    \norm{\epsilon \psi' \circ g(x[\zd])}_{\xLone} &\leq 2\norm{\epsilon \psi' \circ g(x[\ze])}_{\xLone}\leq 2K_g\label{eq:upper_bound_state_mixed}
\end{align}
Now, one can exhibit an upper-bound on the difference $J_\epsilon(x[\zd], \vd)-J_\epsilon(x[\ze], \ue)$ as follows
\begin{equation*}
    J_\epsilon(x[\zd], \vd)-J_\epsilon(x[\ze], \ue) = \Delta_1  + \Delta_2 + \epsilon \Delta_3\label{eq:def_split_c_epsilon}
\end{equation*}
where
\begin{align}
    \Delta_1 &:=
        \varphi(x[\zd](T))-\varphi(x[\ze](T))  + \int_0^T\left[\ell(x[\zd](t),\vd(t))-\ell(x[\ze](t),\ue(t))\right]\xdt\nonumber\\
    \Delta_2 & := \epsilon\int_0^T \sum_i\left[\psi \circ g_i(x[\zd](t))-\psi \circ g_i(x[\ze](t)) \right]\xdt\nonumber\\
     \Delta_3 &:= \int_0^T\left[\psi\circ c(x[\zd](t), \vd(t))-\psi \circ c(x[\ze](t), \ue(t))\right]\xdt\label{eq:def_delta_3_c}
\end{align}
Now, let us upper-bound $\Delta_1$
\begin{multline}
    \Delta_1 \leq 
    \int_0^T\const(\ell)\left(\norm{x[z](t)-x[\ze](t)}+\norm{v(t)-\ue(t)}\right)\xdt
    +\const(\varphi)\norm{\Delta x}_{\xLinfty} \leq \const(\ell,f,g, \varphi, T, R_v, R_x)\label{eq:upperbound_delta1_c}
\end{multline}
Now, let us upper-bound $\Delta_2$.
\begin{equation*}
    \Delta_2 = \epsilon\int_0^T \sum_i\left[\psi \circ g_i(x[\zd](t))-\psi \circ g_i(x[\ze](t)) \right]\xdt = \epsilon\sum_i \int_0^T\int_0^1\psi'\left[g_i(x[\ze](t)) + s \Delta g_i(t)\right]\Delta g_i(t) \xd s \xdt 
\end{equation*}
From the mean value theorem, \cref{eq:interiority_state_mixed,eq:upper_bound_state_mixed}, $\exists \theta_t\in[0,1]$ such that
\begin{equation}
    \Delta_2 = \epsilon\sum_i \int_0^T\psi'\circ g_i(x[\ze](t) + \theta_t \Delta g(t))\Delta g(t) \xdt \leq \sum_i 2 K_g \norm{\Delta g}_{\xLinfty}\leq \const(f,g,K_g, R_x,T)\label{eq:upperbound_delta2_c}
\end{equation}
Now, let us upper-bound $\Delta_3$ defined in \cref{eq:def_delta_3_c}. To do so, let us introduce the following useful subsets of $[0,T]$
 \begin{align*}
    E_1&:=(c(x[\ze], \ue))^{-1}\left((-\infty,-C_\delta]\right)\\
    E_2&:=(c(x[\ze], \ue))^{-1}\left((-C_\delta,0]\right)
\end{align*}
Let us decompose $\Delta_3$ as follows $\Delta_3 := \Delta_{3,1} + \Delta_{3,2}$, with
\begin{equation*}
    \Delta_{3,i} := \int_{E_i} \psi \circ c(x[\zd](t), \vd(t)) - \psi \circ c(x[\ze](t), \ue(t))\xdt,\;i=1,2
\end{equation*}
By convexity of the $\log$ penalty, i.e. $\psi$, we have
\begin{equation}
    \label{eq:upperbound_delta31}
    \Delta_{3,1}\leq\int_{E_1} \psi'(C_\delta) \norm{\Delta c}_{\xLinfty} \xdt \leq 
    \const(T,f,c,\psi,C_\delta)\norm{\zd - \ze}_{\Vad}
    \leq \const(T, f, c, \psi, C_\delta, R_v )
\end{equation}
In addition, we have 
\begin{equation*}
     \Delta_{3,2} =
     \int_{E_2}\left(\int_{0}^{1}\psi'(c(x[\ze](t), \ue(t)) + s\Delta c(t))\Delta c(t) \xd s\right) \xdt
\end{equation*}
Since $\forall t\in E_2, \,\Delta c(t)  < -C_\delta$, we also have
\begin{equation}
    \label{eq:majoS2_c}
    \Delta_{3,2} \leq -C_\delta\int_{E_2}\left(\int_{0}^{1}\psi'(c(x[\ze](t), \ue(t)) + s\Delta c(t))\xd s\right) \xdt
\end{equation}
From the mean value theorem, $\forall t\in E_2$, $\exists \sigma_t\in[0,1]$ such that
\begin{equation}
\label{eq:mean_value_c}
    \psi'\left(c(x[\ze](t), \ue(t)) + \sigma_t \Delta c(t)\right)= \int_{0}^{1} \psi'(c(x[\ze](t), \ue(t)) + s\Delta c(t))\xd s
\end{equation}
Since for all $t\in E_2$, we have $c(x[\ze)(t), \ue(t)) - c(x[\zd](t), \vd(t))\geq C_\delta$ and since $\psi'$ is strictly increasing we have $\sigma_t \in (0,1)$ and 
\begin{equation*}
    \psi'\circ c(x[\zd](t), \vd(t)) < \psi'(c(x[\ze](t), \ue(t)) + \sigma_t \Delta c(t)) <\psi'\circ c(x[\ze](t), \ue(t))
\end{equation*}
From the intermediate value theorem, $\exists \bar{\sigma}\in(0,1)$ such that $\forall t\in E_2$ we have
\begin{equation}
\label{eq:lower_bound_intermediate_c}
    \psi'(c(x[\ze](t), \ue(t)) + \sigma_t \Delta c(t)) \geq (1 - \bar{\sigma})\psi'\circ c(x[\ze](t), \ue(t)) + \bar{\sigma} \psi'\circ c(x[\zd](t), \vd(t))
\end{equation}
Gathering \cref{eq:majoS2_c,eq:mean_value_c,eq:lower_bound_intermediate_c} yields
\begin{align}
    \Delta_{3,2} &\leq
    - C_\delta \int_{E_2}\bigg( (1-\bar{\sigma})\psi' \circ c(x[\ze](t), \ue(t)) + \bar{\sigma} \psi' \circ c(x[\zd](t), \vd(t)) \bigg)\xdt\nonumber\\
    &\leq - C_\delta\left( (1 - \bar{\sigma})\int_{E_2}\psi' \circ c(x[\ze](t), \ue(t))\xdt + \const(\bar{\sigma}, \psi, T,c, C_\delta)\right)\label{eq:upperbound_delta32}
\end{align}
Gathering \cref{eq:upperbound_delta31,eq:upperbound_delta32} we have 
\begin{equation}
    \label{eq:upperbound_delta3_c}
    \epsilon \Delta_3 \leq \const(T,f,c,C_\delta,\epsilon_0, \psi, \bar{\sigma},R_v) -\epsilon C_\delta(1 - \bar{\sigma})\int_{E_2}\psi' \circ c(x[\ze](t), \ue(t))\xdt
\end{equation}
Gathering \cref{eq:upperbound_delta1_c,eq:upperbound_delta2_c,eq:upperbound_delta3_c} yields
\begin{multline}
    J_\epsilon(\zd) - J_\epsilon(\ze) = \Delta_1 + \Delta_2 + \epsilon \Delta_3 \\ \leq
    \const(\ell,f,g, \varphi, T, K_g,R_v, R_x, C_\delta,c,\epsilon_0, \psi, \bar{\sigma})-\epsilon C_\delta(1 - \bar{\sigma})\int_{E_2}\psi' \circ c(x[\ze](t), \ue(t))\xdt
    \label{eq:upperbound_J2}
\end{multline}
Now let us prove \cref{eq:boundedness_lambda_c} by contradiction and assume that
\begin{equation*}
    \forall K_c>0, \exists \epsilon>0 \st \norm{\epsilon \psi'\circ c(x[\ze], \ue) }_{\xLone}>K_c
\end{equation*}
From the definition of $E_1$ and $E_2$, we have
\begin{equation*}
    \norm{\epsilon \psi'(c(x[\ze], \ue)}_{\xLone} = 
        \int_{E_1}\epsilon\psi'(c(x[\ze](t), \ue(t))\xdt +
        \int_{ E_2}\epsilon \psi'(c(x[\ze](t), \ue(t))\xdt
\end{equation*}
which, in turns yields
\begin{equation}
    \int_{ E_2}\epsilon \psi'(c(x[\ze](t), \ue(t))\xdt > K_c - \epsilon_0\psi'(C_\delta)T
    \label{eq:eq_majo_int_e2_rho_c}
\end{equation}
gathering \cref{eq:upperbound_J2,eq:eq_majo_int_e2_rho_c} yields
\begin{equation*}
    J_\epsilon(\zd) - J_\epsilon(\ze) \leq\\
    \const(\ell,f,g, \varphi, T, K_g,R_v, R_x,C_\delta,c,\epsilon_0, \psi, \bar{\sigma})- C_\delta(1 - \bar{\sigma})K_c
\end{equation*}
For $K_c$ large enough, $J_\epsilon(\zd) - J_\epsilon(\ze)<0$, which contradicts the optimality of $(x[\ze], \ue)$, proves \cref{eq:boundedness_lambda_c} and concludes the proof.
\end{proof}
\subsection{Proof of \cref{eq:boundedness_p}}
\begin{prpstn}
\label{prop:well_posed_p}
There exists a constant $K_p<+\infty$ such that for all $\epsilon \in  (0,\epsilon_0)$ and for any $(x[\ue, x^0_\epsilon], \ue)$ locally optimal solution of Problem \cref{eq:def_log_barrier_ocp}, equation \cref{eq:boundedness_p} holds.
\end{prpstn}

\begin{proof}
First, using \cref{eq:penalized_pontryagin_extremal_2} one has
\begin{multline*}
    \norm{\bpe(T)- \bpe(s)} 
    \leq \int_s^T\left(\norm{\ell_x'(\bxe,\bue)}_{\xLinfty}+\norm{f_x'(\bxe,\bue)}_{\xLinfty}\norm{\bpe(t)}\right)\xdt +\sum_i\norm{g_i'(\bxe)}_{\xLinfty}\norm{\epsilon \psi'\circ g_i(\bxe)}_{\xLone}\\
    +\sum_i\norm{c_{i,x}'(\bxe, \bue)}_{\xLinfty}\norm{\epsilon \psi'\circ c_i(\bxe,\bue)}_{\xLone}
\end{multline*}
From the continuity of $\ell_x'$, $f_x'$, $g_i'$, $c'_{i,x}$ and since $\bxe$ and $\bue$ are bounded, we have $\norm{\ell_x'(\bxe,\bue)}_{\xLinfty}<\const(\ell,f)$ and  $\norm{f_x'(\bxe,\bue)}_{\xLinfty}<\const(f)$. In addition, the terms on the right-hand side of \cref{eq:penalized_pontryagin_extremal_6} are bounded which yields that $\norm{\bpe(T)}\leq\const(f,h)$. The derivatives of the penalty functions being uniformly $\xLone$-bounded one can use Grönwall Lemma which proves that $\exists K_p<+\infty$ such that $\forall s\in[0,T]$ we have $\norm{\bpe(s)}\leq K_p$ which concludes the proof.
\end{proof}
\section{Proof of \cref{lem:conv_penalties}}
\label{sec:add_proof}

\subsection{Proof of $i)$}
\label{sec:proof_conv_mu}
From \cref{thm:interior_all}, the sequence $(\lambda^g_{\epsilon_n})_n$ is uniformly $\xLone$-bounded. This sequence can be canonically injected into $\mathcal{M}([0,T])^{n_g}$ and is still uniformly bounded with respect to $\mathcal{M}([0,T])^{n_g}$ strong topology. From the weak $\ast$ compactness of the unit ball of $\mathcal{M}([0,T])^{n_g}$, there exists a subsequence and a measure $\xd \bar{\mu}\in \mathcal{M}([0,T])^{n_g}$ such that 
$$\lim_{k\rightarrow +\infty} \lambda^g_{\epsilon_{n_k}}\xd t\stackrel{\ast}{\rightharpoonup} \xd \bar{\mu}$$

\subsection{Proof of $ii)$.}
\label{sec:proof_of_2}
For all $K_c>0$, assume that there exists $E\subseteq[0,T]$ of strictly positive measure such that $I^c_{\bue,\bxe(0)}(t,K_c/\epsilon )\neq\emptyset$ for all $t\in E$. Now, let us denote $C(t):= c_{I^c_{\bue,\bxe(0)}(t,K_c/\epsilon ),u}'(\bxe, \bue)$ and let us define $v$ as follows
    \begin{equation*}
        v(t) := \begin{cases}
            - C(t)^\top\left[C(t) C(t)^\top\right]^{-1}.e_v,\forall t\in E\\
            0\textrm{ otherwise}
        \end{cases}
    \end{equation*}
   where $\R^{\vert I^c_{\bue,\bxe(0)}(t,K_c/\epsilon )\vert} \ni e_v :=\begin{pmatrix}\epsilon/K_c &\dots & \epsilon/K_c\end{pmatrix}^\top$. Since $\bue$ is a locally optimal solution of Problem \cref{eq:def_log_barrier_ocp}, and using \cref{def:penalized_pre_hamiltonian}, we have for almost all $t\in E$
    \begin{multline*}
        H^\psi(\bxe(t), \bue(t)+ v(t), \bpe(t), \epsilon) - H^\psi(\bxe(t), \bue(t), \bpe(t), \epsilon) = H(\bxe(t), \bue(t) + v(t) ,\bpe(t)) - H(\bxe(t), \bue(t),\bpe(t))\\
        + \epsilon\sum_{i\in I^c_{\bue,\bxe(0)}(t,K_c/\epsilon)} \log\left(\frac{c_i(\bxe(t), \bue(t))}{c_i(\bxe(t), \bue(t) + v(t))}\right)
    \end{multline*}
    From the mean value theorem, $\exists s\in[0,1]$ such that 
    \begin{multline*}
         H^\psi(\bxe(t), \bue(t)+ v(t), \bpe(t), \epsilon) - H^\psi(\bxe(t), \bue(t) , \bpe(t), \epsilon) \leq \\
        H'_{u}(\bxe(t), \bue(t) + s v(t),\bpe(t)). v(t)  -\epsilon\sum_{i\in I^c_{\bue,\bxe(0)}(t,K_c/\epsilon)} \log(2)
    \end{multline*}
    Since $\ell, f$ are at least $\xCtwo$ and since $C(t)$ is $\xLinfty$-bounded by some constant $M$, from \cref{eq:boundedness_p} we have
    \begin{multline*}
     H^\psi(\bxe(t), \bue(t)+ v(t), \bpe(t), \epsilon) - H^\psi(\bxe(t), \bue(t) , \bpe(t), \epsilon) \leq\\
     \epsilon\left(\frac{\const(\ell, f, T,R_v, R_x, K_p,M)}{K_c}-\sum_{i\in I^c_{\bue,\bxe(0)}(t,K_c/\epsilon)} \log(2) \right)
    \end{multline*}
    which is negative for $K_c$ large enough, contradicts the local optimality of $\bue$ and proves the result.

\subsection{Proof of $iii)$.}
\label{sec:proof_of_3}
The proof of this result consists in proving that the mapping 
\begin{equation*}
    \Lambda_{\epsilon}^{c}: \xLone([0,T];\R^{n_c}) \ni w \mapsto \int_0^T \lambda_{\epsilon}^{c}(t).w(t)\xdt \in \R
\end{equation*}
is a continuous linear form on $\xLone([0,T];\R^{n_c})$. From \cref{eq:interiority_mixed}, $c(\bxe, \bue)$ is strictly negative (not active) almost everywhere. Therefore the Hamiltonian minimization condition of the Pontryagin maximum principle writes ${H^\psi_u}'(\bxe,\bue, \bpe, \epsilon) = 0$ for almost all time. Then, for all $ v\in \xLone([0,T];\R^m)$, one has
    \begin{equation}
        \left\vert \int_0^T \sum_{i=1}^{n_c} \lambda_{\epsilon}^{c_i}(t) c'_{i, u}(\bxe(t), \bue(t)).v(t) \xdt \right\vert \leq 
        \norm{ H_u'(\bxe,\bue, \bpe)}_{\xLinfty}\norm{v}_{\xLone} \leq \const(f,\ell,h)\norm{v}_{\xLone}
        \label{eq:Linfty_L1_pen_mixed}
    \end{equation}
    Now, let us denote $C(t):= c_{I^c_{\bue, \bx^0_{\epsilon}}(t,n),u}'(\bxe, \bue)$ and for all $w\in \xLinfty([0,T];\R^{n_c})$, let us define $v\in \xLinfty([0,T];\R^m)$ as follows
    \begin{equation}
        v(t) :=\begin{cases}
                C(t)^\top\left[C(t) C(t)^\top\right]^{-1}w_{I^c_{\bue, \bx^0_ {\epsilon}}(t,n)}(t) &\textrm{if } I^c_{\bue, \bx^0_ {\epsilon}}(t,n)\neq \emptyset\\
                0 &\textrm{otherwise}
        \end{cases}
        \label{eq:def_vt_from_wt}
    \end{equation}
    Since $C(t)$ is $\xLinfty$-bounded, there exists $M\in (0, +\infty)$ such that $\norm{v}_{\xLone}\leq M \norm{w}_{\xLone}$ 
    In addition, let us define $\theta_\epsilon\in \xLone([0,T];\R^{n_c})$ as follows
    \begin{equation}
        \theta_{\epsilon,i}(t):= \begin{cases}
            \lambda_{\epsilon}^{c_i}(t) &\textrm{if } i\in I^c_{\bue, \bx^0_{\epsilon}}(t,n)\\
            0 &\textrm{otherwise}
        \end{cases}
        \label{eq:def_lambda_epsilon_t}
    \end{equation}
    Gathering \cref{eq:Linfty_L1_pen_mixed,eq:def_vt_from_wt,eq:def_lambda_epsilon_t} yields
     \begin{equation}
       \left\vert  \int_0^T \sum_{i=1}^{n_c}\lambda_{\epsilon}^{c_i}(t) c'_{i,u}(\bxe(t), \bue(t)).v(t) \xdt \right\vert = \left\vert \int_0^T \theta_\epsilon(t).w(t) \xdt \right\vert
       \label{eq:left_right_equal_lambda_w_mixed}
    \end{equation}
    Gathering \cref{eq:Linfty_L1_pen_mixed,eq:left_right_equal_lambda_w_mixed} we have for all $w\in \xLinfty([0,T];\R^{n_c})$
    \begin{equation*}
        \left\vert \int_0^T \theta_\epsilon(t).w(t) \xdt \right\vert =\left\vert  \int_0^T \sum_{i=1}^{n_c}\lambda_{\epsilon}^{c_i}(t) c'_{i,u}(\bxe(t), \bue(t)).v(t) \xdt \right\vert \leq \const(f,\ell,h)\norm{v}_{\xLone}\leq \const(f,\ell,h)\norm{w}_{\xLone}
    \end{equation*}
   Since $\xLinfty([0,T];\R^{n_c})$ is dense in $\xLone([0,T];\R^{n_c})$, $\theta_\epsilon$ can be extended to a continuous linear form over $\xLone([0,T];\R^{n_c})$. Therefore $\theta_\epsilon \in \xLone([0,T];\R^{n_c})^*\cong \xLinfty([0,T];\R^{n_c})$, hence $\lambda_\epsilon^{c_i} \in \xLinfty([0,T];\R_+)$ which concludes the proof.

\bibliographystyle{plain}
\bibliography{references}

\begin{thebibliography}{10}

\bibitem{Bogachev.2007}
V.I. Bogachev.
\newblock {\em Measure Theory}.
\newblock Springer Berlin Heidelberg, Berlin, Heidelberg, 2007.

\bibitem{Bonnans.2014}
J.~F. Bonnans, X.~Dupuis, and L.~Pfeiffer.
\newblock Second-order necessary conditions in pontryagin form for optimal
  control problems.
\newblock {\em SIAM Journal on Control and Optimization}, 52(6):3887--3916,
  2014.

\bibitem{FredericBonnans.2014}
J.F. {Bonnans}, Xavier Dupuis, and Laurent Pfeiffer.
\newblock Second-order sufficient conditions for strong solutions to optimal
  control problems.
\newblock {\em ESAIM: Control, Optimisation and Calculus of Variations},
  20(3):704--724, 2014.

\bibitem{Bonnans_log}
J.F. Bonnans and T.~Guilbaud.
\newblock Using logarithmic penalties in the shooting algorithm for optimal
  control problems.
\newblock {\em Optimal Control Applications and Methods}, 24:257--278, 2003.

\bibitem{BonnansHermant}
J.F. Bonnans and A.~Hermant.
\newblock Second-order analysis for optimal control problems with pure state
  constraints and mixed control-state constraints.
\newblock {\em Ann. Inst. H. Poincaré Anal. Non Linéaire}, 26(2):561–598,
  2009.

\bibitem{Bonnans.Shapiro}
J.F. Bonnans and A.~Shapiro.
\newblock {\em Perturbation analysis of optimization problems}.
\newblock Springer-Verlag, New York, 2000.

\bibitem{brezis}
H.~Brézis.
\newblock {\em Functional Analysis, Sobolev Spaces and Partial Differential
  Equations}.
\newblock Springer, 2010.

\bibitem{CAILLAU2023559}
J-B. Caillau, R.~Ferretti, E.~Trélat, and H.~Zidani.
\newblock Chapter 15 - an algorithmic guide for finite-dimensional optimal
  control problems.
\newblock In {\em Numerical Control: Part B}, volume~24 of {\em Handbook of
  Numerical Analysis}, pages 559--626. Elsevier, 2023.

\bibitem{graichen}
K.~Graichen and N.~Petit.
\newblock Incorporating a class of constraints into the dynamics of optimal
  control problems.
\newblock {\em Optimal Control Applications and Methods}, 30:537--561, 2009.

\bibitem{Hartl}
R.~F. Hartl, Suresh Sethi, and Raymond Vickson.
\newblock A survey of the maximum principles for optimal control problems with
  state constraints.
\newblock {\em SIAM Review}, 37(2):181--218, 1995.

\bibitem{Khalil}
H.~Khalil.
\newblock {\em Non {L}inear {S}ystems}.
\newblock Prentice {H}all, 2002.

\bibitem{Kierzenka2001ABS}
J.~Kierzenka and L.F. Shampine.
\newblock A bvp solver based on residual control and the maltab pse.
\newblock {\em ACM Trans. Math. Softw.}, 27:299--316, 2001.

\bibitem{kurdila}
A.J. Kurdila and M.~Zabarankin.
\newblock {\em Convex Functional Analysis}.
\newblock {Birkh{\"a}user Boston}, 2005.

\bibitem{lasdon67}
L.~Lasdon, A.~Waren, and R.~Rice.
\newblock An interior penalty method for inequality constrained optimal control
  problems.
\newblock {\em IEEE Transactions on Automatic Control}, 12:388--395, 1967.

\bibitem{maliSIAM}
P.~Malisani.
\newblock Interior point methods in optimal control problems for affine
  systems: Convergence results and solving algorithms.
\newblock {\em SIAM Journal on Control and Optimization}, 61(6), 2023.

\bibitem{source_ifpen}
P.~Malisani.
\newblock Python source code for "{I}nterior point methods in optimal control".
\newblock \url{https://ifpen-gitlab.appcollaboratif.fr/detocs/ipm_ocp}, 2023.

\bibitem{maliOCAM}
P.~Malisani, F.~Chaplais, and N.~Petit.
\newblock An interior penalty method for optimal control problems with state
  and input constraints of nonlinear systems.
\newblock {\em Optimal Control Applications and Methods}, 37:3--33, 2014.

\bibitem{Maurer}
A.~Maurer and J.~Zowe.
\newblock First and second-order necessary and sufficient optimality conditions
  for infinite-dimensional programming problems.
\newblock {\em Math. Programming}, 16:98--110, 1979.

\bibitem{Renardy.2004}
M.~Renardy and R.~Rogers.
\newblock {\em An introduction to partial differential equations}, volume~13 of
  {\em Texts in applied mathematics}.
\newblock Springer, New York, second edition edition, 2004.

\bibitem{Rudin}
W.~Rudin.
\newblock {\em Principles of Mathematical Analysis}.
\newblock McGraw-Hill, 3rd edition, 1976.

\bibitem{seywald1993goddard}
H.~Seywald and E.M. Cliff.
\newblock Goddard problem in presence of a dynamic pressure limit.
\newblock {\em Journal of Guidance, Control, and Dynamics}, 16(4):776--781,
  1993.

\bibitem{Weiser}
M.~Weiser.
\newblock Interior point methods in function space.
\newblock {\em SIAM Journal on Control and Optimization}, 44(5):1766--1786,
  2005.

\end{thebibliography}

\end{document}